\documentclass[12pt,twoside]{article}
\usepackage{amsmath,amssymb}
\usepackage{graphicx}
\usepackage{color}

\makeatletter 
\@addtoreset{equation}{section}
\makeatother  

\newtheorem{df}{Definition}[section]
\newtheorem{lm}{Lemma}[section]
\newtheorem{thm}{Theorem}[section]
\newtheorem{rmk}{Remark}[section]

\newtheorem{pro}{Proposition}[section]
\newtheorem{co}{Corollary}[section]

\makeatletter
\evensidemargin
\oddsidemargin
\makeatother

\title{\Large\bf Mittag-Leffler Dichotomy and Roughness of Conformable Fractional Differential Equations
\thanks{Supported by NSFC \#11701400 and \#12071317.}
}

\author{{\sc Baishun Wang},~~~{\sc Jun Zhou}
\footnote{E-mail address:matzhj@163.com.}
\\
{\small School of Mathematical Science, Sichuan Normal University,}\\
{\small Chengdu, Sichuan 610066, P. R. China}}

\date{}

\begin{document}
\thispagestyle{empty}
\setcounter{page}{1}

\maketitle

\begin{abstract}

The solutions of traditional fractional differential equations neither satisfy group property nor generate dynamical systems,
so the study on hyperbolicity is in blank.
Relying on the new proposed conformable fractional derivative,
we investigate dichotomy of conformable fractional equations,
including structure of solutions of linear systems, Mittag-Leffler dichotomy and stability,
roughness and nonuniform dichotomy.

\\
\vskip 0.2cm
{\bf AMS (2020) subject classification:} 34D09; 34A08; 37D25

\end{abstract}


\section{Introduction}
\setcounter{equation}{0}

The well-known dichotomy concept on various of hyperbolic systems, e.g. ODEs
\begin{eqnarray}
\label{xAx}
\dot{x}=A(t)x,~~~t\in J,
\end{eqnarray}
where interval $J\subset\mathbb{R}$, is said that there exist a projection
matrix $P$ and a fundamental matrix $X(t)$ of \eqref{xAx}, and positive constants
$K_i$ and $\beta_i$ ($i=1,2$) such that for all $t,s\in J$,
\begin{eqnarray}
\begin{split}
\|X(t)PX^{-1}(s)\|\leq K_1e^{-\beta_1(t-s)},~~~t\geq s,\\
\|X(t)(I-P)X^{-1}(s)\|\leq K_2e^{-\beta_2(s-t)},~~~s\geq t.
\end{split}
\label{XPXIP}
\end{eqnarray}
Correspondingly, the roughness of dichotomy is regarded as the persistence of dichotomy undergoing
small linear perturbation, i.e., the perturbed system
\begin{eqnarray*}
\dot{x}=(A(t)+B(t))x,~~~t\in J
\end{eqnarray*}
still admits dichotomy behaviour in the form of \eqref{XPXIP}
along with small variations of $P$, $X(t)$, $K_i$ and $\beta_i$ ($i=1,2$).
According to the difference of asymptotic rate,
there are diverse dichotomies, e.g.,
the classical exponential dichotomy \eqref{XPXIP} (\cite{C78}),
$(h,k)$-dichotomy (\cite{NP95}), polynomial dichotomy (\cite{BV09}), etc.
The dichotomies and their corresponding properties are core issues
in the field of dynamical systems,
which can be traced back to the papers of Perron (\cite{P30}) and Li (\cite{L34})
on conditional stability of linear differential and difference equations respectively.
And they are gradually formalized, developed and summarized in literatures \cite{MS58,MS60,C78}.
Recently decades, many research results were devoted to exploring the existence criteria of exponential dichotomy
(see Hale (\cite{H69}), Chow and Leiva (\cite{CL95}), Sasu (\cite{S08}),
Barreira and Valls (\cite{BV08}), Battelli and Palmer (\cite{BP19}) and the references therein).
The roughness referred before are also widely focused on,
and firstly demonstrated by Massera and Sch\"affer (\cite{MS58})
under hypothesis that the original matrix $A(t)$ is bounded.
Sch\"affer (\cite{S63}) subsequently eliminated the assumption of boundedness.
Coppel (\cite{C67}) gave a general elementary proof of roughness if matrix $A(t)$ commutes
with the projection $P$.
In 1978 Coppel (\cite[pp.28-33]{C78}) exhibited a simpler proof via the so-called {\it projected integral inequalities}
raised by Hale (\cite[pp.110-111]{H69}).
Later, Naulin and Pinto (\cite{NP98}) improved the size of perturbation $B(t)$ in Coppel's \cite[pp.34-35]{C78}
without boundedness of $A(t)$ yet.
Popescu (\cite{P06}) further generalized the results of \cite{C78} and \cite{NP98}
to infinite dimensional Banach spaces.
Thereafter, the notion of nonuniform exponential dichotomy,
roughly speaking dichotomy formula \eqref{XPXIP} involving extra nonuniform constants
in exponents, was proposed by Barreira and Valls (\cite{BV08}),
where its roughness was also studied.
In 2013 Zhou, Lu and Zhang (\cite{ZLZ13}) studied the roughness of tempered exponential dichotomy
for random difference equations in Banach spaces lack of the so-called {\it Multiplicative Ergodic Theorem}.
Moreover, plenty of works on the roughness of exponential dichotomy could be found
in \cite{CL99,CL95,C78,NP98,P06,BV08} for continuous dynamical systems and
in \cite{PS99,S08,ZLZ13,ZZ16} for discrete dynamical systems and references therein.
In addition, the corresponding admissibility problem of dichotomy,
i.e. admissible functions pair of solutions $x$ and inhomogeneous perturbations $f$,
was investigated extensively in \cite{MS60,S08,ZZ16,BDV18,DZZ22} and so on.

Although the research on dichotomy involved ODEs (\cite{MS58,MS60,H69,C78,NP95,NP98,CL99,P06,BV08,BV09,BP19,DZZ22}),
difference equations (\cite{L34,S08,ZZ16,ZLZ13}), functional differential equations (\cite{P71}),
random systems (\cite{A98,G95,C97,LL10}), skew-product semiflows (\cite{CL95,PS99}), etc.,
till now there is no result of dichotomy for fractional differential equations (FDEs for short).
Fractional derivative started from a letter from L'Hospital to Leibniz about discussing the meaning of a half derivative.
From then on, because of better approximation to practical model associated with memory and hereditary phenomena
than ODEs and PDEs, FDEs are steadily developed in the aspects of Physics and Chemistry (\cite{ZCS20,We07}),
Biology and Medicine (\cite{Di13,SSM18}), Engineering and Control Theory (\cite{BF07,P99})
and Economics and Psychology (\cite{CC07,SXY10}), especially in recent decades (see monographs \cite{P99,KST06,GPH11,ZWZ16,J21,DMD19}).
Traditional definitions of fractional derivative and integral, such as Riemann-Liouville's,
Caputo's and Gr\"unwald-Letnikov's (\cite{P99,KST06}), have no product rule and chain rule of derivative,
such that the solutions of FDEs neither fulfil group property nor generate dynamical systems.
There is a vast of works on well-posedness (\cite{ZWZ16}), stability (\cite{KST06}),
Laplace transform method and optimal control (\cite{P99}), variational method, attractors and
numerical solutions (\cite{GPH11}) and chaos (\cite{BWZ17}) of FDEs,
but the study on hyperbolicity of FDEs is temporarily in blank state.
Until 2014, Khalil et al. (\cite{KHY14}) introduced a new definition of fractional derivative,
that is so-called {\it conformable fractional derivative},
which can almost satisfy all corresponding characteristics analogous to integer derivative.
Thus, the solutions of CFDEs (abbreviation of conformable fractional differential equations) also can generate dynamical systems,
which makes it possible to consider the hyperbolic behaviors of FDEs.
Later, Abdeljawad (\cite{A15}) accomplished the definition of left and right conformable fractional derivatives
and the variation of constants formula of CFDEs
and solved CFDEs via Laplace transform.
In 2017 Souahi et al. (\cite{SMH17}) employed Lyapunov direct method to present the stability,
asymptotic stability and exponential stability of CFDEs.
In 2019 Khan et al. (\cite{UA19}) further verified
the generalized definition and its semigroup and linear properties of conformable derivative and
existence and uniqueness of solutions for CFDEs.
During the same year, Balci et al. (\cite{BOK19}) displayed the Neimark-Sacker bifurcation and chaotic behavior
for a tumor-immune system modelled by a CFDE.
In 2020 Xie et al. (\cite{XLW20}) showed an exact solution and difference scheme for a gray model with conformable derivative.
Recently, Wu et al. (\cite{WZL22}) revealed the Hyers-Ulam stability of a conformable fractional model.

In this paper we attempt to establish the theory of dichotomy for CFDEs.
In order to generalize the hyperbolicity of ODEs to CFDEs,
we first modify the definitions of conformable fractional derivative and integral
and a Mittag-Leffler-type function originated from \cite{KHY14,A15}.
Subsequently, we derive the well-posedness of solutions,
conformable integral inequality and variation of constants formula for CFDEs
and structure of solutions and operator semigroups for linear CFDEs.
These fundamental theories are achieved in section 2.
In section 3 we provide the definitions of so-called {\it Mittag-Leffler stability and dichotomy}
with respect to CFDEs, whose asymptotic rate is a Mittag-Leffler-type function.
These stability and dichotomy include the classical exponential stability and dichotomy
(\cite{C78}) in ODEs with integer derivative as special cases.
Meanwhile, we develop the conformable fractional integral versions
of projected inequalities to prove the existence of Mittag-Leffler dichotomy
and corresponding invariant manifolds.
In section 4 we discuss the roughness of Mittag-Leffler dichotomy in $\mathbb{R}_+$.
In section 5 we additionally study nonuniform Mittag-Leffler stability and dichotomy and their roughness in $\mathbb{R}_+$ too.
Our results extend the works of Hale (\cite{H69}), Coppel (\cite{C78}), Barreira and Valls (\cite{BV08}) to CFDEs.

\section{Linear CFDEs}
\setcounter{equation}{0}

Throughout this paper, we define the following functions sets:
\begin{eqnarray*}
I(\mathbb{R},\mathbb{R})
    :=\!\!\!\!\!\!\!\!&& \{\varphi:\mathbb{R}\to \mathbb{R}\,|\,
    \varphi \mbox{ is a nondecreasing function}\},
    \\
C(\mathbb{R},\mathbb{R})
    :=\!\!\!\!\!\!\!\!&& \{\varphi:\mathbb{R}\to \mathbb{R}\,|\,
    \varphi \mbox{ is a continuous function}\},
    \\
C_b(\mathbb{R},\mathbb{R})
    :=\!\!\!\!\!\!\!\!&& \{\varphi:\mathbb{R}\to \mathbb{R}\,|\,
    \varphi \mbox{ is a continuous and bounded function}\},
\\
C_I(\mathbb{R},\mathbb{R})
    :=\!\!\!\!\!\!\!\!&& \{\varphi:\mathbb{R}\to \mathbb{R}\,|\,
    \varphi \mbox{ is a continuous and nondecreasing function}\}.
\end{eqnarray*}
Further, set constants $\alpha\in(0,1]$ and $t_*,t_0,t^*$ satisfying $t_*<t_0<t^*$,
function $f:(t_*,t^*)\to\mathbb{R}$, and norms
\begin{eqnarray*}
\|x(t)\|\!\!\!\!\!\!\!\!&&:=\sum\limits^n_{i=1}|x_i(t)|,~~~x:[t_0,+\infty)\to\mathbb{R}^n,
\\
\|A(t)\|\!\!\!\!\!\!\!\!&&:=\max\{\sum\limits^n_{i=1}|a_{i1}(t)|,\sum\limits^n_{i=1}|a_{i2}(t)|,...,\sum\limits^n_{i=1}|a_{in}(t)|\},
~~~A:\mathbb{R}\to\mathbb{R}^{n\times n}.
\end{eqnarray*}
In this section, we focus on the qualitative properties of linear
CFDEs and their perturbations.
Analogously to the linear ODEs, there also exists fundamental solutions
for linear CFDEs.
Consider the nonautonomous linear CFDE
\begin{eqnarray}
{\cal T}^\alpha x=A(t)x,~~~(t,x)\in\mathbb{R}^{n+1},
\label{TxAx}
\end{eqnarray}
where matrix function $A\in C(\mathbb{R},\mathbb{R}^{n\times n})$.
Our primary purpose is to establish its well-posedness, e.g. existence and uniqueness,
continuous dependence on initial data of solutions and continuation of solutions.
Before this, as a preliminary, we modify the definition and
some properties of conformable fractional derivative and integral
raised by Khalil, Horani, Yousef and Sababheh(\cite{KHY14})
and Abdeljawad(\cite{A15}), to make them make more sense.

\begin{df}
\label{CFD}
The $\alpha$-conformable fractional derivative of $f$ is defined as
\begin{eqnarray}
{\cal T}^\alpha f(t):=\lim\limits_{\varepsilon \to 0}\frac{f(t+\varepsilon |t|^{1-\alpha})-f(t)}{\varepsilon},
~~~t\in(t_*,t^*).
\label{Taf}
\end{eqnarray}
In particular, if $\lim\limits_{t\to 0}{\cal T}^\alpha f(t)$ exists, then
$${\cal T}^\alpha f(0):=\lim_{t\to 0}{\cal T}^\alpha f(t).$$
Here function $f$ is called as $\alpha$-conformable differentiable,
if ${\cal T}^\alpha f(t)$ exists.
\end{df}

Our Definition \ref{CFD} extends the one in \cite{KHY14} to the case of $t\leq0$.
And different from the definition in \cite{A15},
there are same formulae in \eqref{Taf} for both $t\leq0$ and $t\geq0$.
Further, we can deduce the following relations between conformable fractional derivative
and Newton-Leibniz derivative and between conformable fractional integral
and Riemann integral.

\begin{pro}
\label{TaD}
The $\alpha$-conformable fractional derivative of $f$
can be represented as
\begin{eqnarray*}
\label{Taftf}
{\cal T}^\alpha f(t)=|t|^{1-\alpha}f'(t),~~~t\in(t_*,t^*).
\end{eqnarray*}
\end{pro}

This proposition implies that there are same sign between
conformable fractional derivative and Newton-Leibniz derivative.

\begin{pro}
The $\alpha$-conformable fractional integral of $f$ is given by
\begin{eqnarray*}
{\cal I}^\alpha_{t_0}f(t)
=\int_{t_0}^t|s|^{\alpha-1}f(s)ds,
~~~t\in(t_*,t^*).
\end{eqnarray*}
\label{cfdejfgs}
\end{pro}

Next, we present some properties of conformable fractional derivative and integral.

\begin{pro}
Let real functions $f$ and $g$ be $\alpha$-conformable differentiable in $(t_*,t^*)$,
then the following properties hold:
\begin{description}
  \item[(A1)] ${\cal T}^\alpha {\cal I}^\alpha_{t_0}f(t)=f(t)$, ${\cal I}^\alpha_{t_0}{\cal T}^\alpha f(t)=f(t)-f(t_0)$;
  \item[(A2)] ${\cal T}^\alpha(af+bg)=a{\cal T}^\alpha f+b{\cal T}^\alpha g$, for all $a,b\in \mathbb{R}$;
  \item[(A3)] ${\cal T}^\alpha(fg)=f{\cal T}^\alpha g+g{\cal T}^\alpha f$,
  ${\cal T}^\alpha\big{(}\cfrac{f}{g}\big{)}=\cfrac{g{\cal T}^\alpha f-f{\cal T}^\alpha g}{g^2}$;
  \item[(A4)] ${\cal T}^\alpha f\circ g=|g|^{\alpha-1}{\cal T}^\alpha f(g){\cal T}^\alpha g$, if $g:(t_*,t^*)\to(t_*,t^*)$.
\end{description}
\end{pro}
Specially, the conformable fractional derivatives of some elementary functions
are stated as follows:
\begin{eqnarray*}
\begin{array}{lll}
\!\!\!\!\!\!\!\!&(1)~{\cal T}^\alpha 1=0,~~&(2)~{\cal T}^\alpha e^{ct}=c|t|^{1-\alpha}e^{ct},~\forall c\in\mathbb{R},
\\
\!\!\!\!\!\!\!\!&(3)~{\cal T}^\alpha t^p=pt^{p-\alpha},~\forall t\geq0,p\in\mathbb{R}~~&(4)~{\cal T}^\alpha \cfrac{t^\alpha}{\alpha}=1,~\forall t\geq0,
\\
\!\!\!\!\!\!\!\!&(5)~{\cal T}^\alpha \sin(bt)=b|t|^{1-\alpha}\cos(bt),~\forall b\in\mathbb{R},~~&(6)~{\cal T}^\alpha \cos(bt)=-b|t|^{1-\alpha}\sin(bt),~\forall b\in\mathbb{R}.
\end{array}
\end{eqnarray*}

\begin{pro}
\label{Contin}
If function $f(t)$ is $\alpha$-conformable differentiable at $t=\hat{t}$,
then $f$ is continuous at $\hat{t}$.
\end{pro}

{\bf Proof.}
For any $\hat{t}\in (t_*,t^*)$,
since ${\cal T}^\alpha f(\hat{t})$ exists,
we have
\begin{eqnarray*}
\lim\limits_{\varepsilon\to 0}[f(\hat{t}+\varepsilon |\hat{t}|^{1-\alpha})-f(\hat{t})]
\!\!\!\!\!\!\!\!&&=\lim\limits_{\varepsilon\to 0}\frac{f(\hat{t}+\varepsilon |\hat{t}|^{1-\alpha})-f(\hat{t})}{\varepsilon}
\cdot\varepsilon
\\
\!\!\!\!\!\!\!\!&&={\cal T}^\alpha f(\hat{t})
\lim\limits_{\varepsilon\to 0}\varepsilon
=0,
\end{eqnarray*}
where $\lim\limits_{\varepsilon\to 0}\varepsilon |\hat{t}|^{1-\alpha}=0$.
It leads to the continuity of $f$ at $\hat{t}$,
and the proposition is proved.
$\qquad\Box$

The following special function and fractional integral inequality
both will be useful throughout this paper.

\begin{df}
The following special function is called as a Mittag-Leffler-type function:
\begin{eqnarray*}
E_\alpha (\lambda,t):=
\left\{
\begin{array}{ll}
\exp\Big{(}\lambda \cfrac{t^\alpha}{\alpha}\Big{)}
=\sum\limits_{k=0}^{+\infty}\cfrac{\lambda ^k t^{\alpha k}}{\alpha ^kk!},
~~~\lambda \in \mathbb{R},~~~t\in\mathbb{R}_+,
\\
\exp\Big{(}-\lambda \cfrac{(-t)^\alpha}{\alpha}\Big{)}
=\sum\limits_{k=0}^{+\infty}\cfrac{(-\lambda)^k(-t)^{\alpha k}}{\alpha ^kk!},
~~~\lambda \in \mathbb{R},~~~t\in\mathbb{R}_-.
\end{array}
\right.
\end{eqnarray*}
\label{MLfunc}
\end{df}

\begin{lm}
Let functions $a\in I([t_0,t^*),\mathbb{R}_+)$ and
$f\in C([t_0,t^*),\mathbb{R}_+)$.
Assume that $u:[t_0,t^*)\to\mathbb{R}_+$ satisfies fractional
integral inequality
\begin{eqnarray}
u(t)\leq a(t)+{\cal I}^\alpha_{t_0}f(t)u(t),~~~t\in[t_0,t^*).
\label{uaIf}
\end{eqnarray}
Then $u$ can be estimated by
\begin{eqnarray}
u(t)\!\!\!\!&\leq\!\!\!\!& a(t)e^{{\cal I}^\alpha_{t_0}f(t)}
\nonumber
\\
\!\!\!\!&\leq\!\!\!\!& a(t)E_\alpha\Big{(}\sup\limits_{s\in[t_0,t]}f(s),|t|\Big{)}E_\alpha\Big{(}\sup\limits_{s\in[t_0,t]}f(s),|t_0|\Big{)},
~~~t\in[t_0,t^*).
\label{uaeIf}
\end{eqnarray}
\label{XEE}
\end{lm}

{\bf Proof.}
Dividing $a(t)$ on the both hands sides of \eqref{uaIf},
by the monotonicity of $a$, we derive
\begin{eqnarray}
\frac{u(t)}{a(t)}\leq 1+{\cal I}^\alpha_{t_0}f(t)\frac{u(t)}{a(t)},
~~~t\in[t_0,t^*).
\label{ua1sf}
\end{eqnarray}
Let $v(t):=1+{\cal I}^\alpha_{t_0}f(t)\cfrac{u(t)}{a(t)}$,
then
\begin{eqnarray*}
v'(t)\leq |t|^{\alpha-1}f(t)v(t),~~~t\in[t_0,t^*).
\end{eqnarray*}
Integrating both hands sides of equation above from $t_0$ to $t$,
we get
\begin{eqnarray*}
\ln{v(t)}\leq {\cal I}^\alpha_{t_0}f(t),~~~t\in[t_0,t^*),
\end{eqnarray*}
which implies
\begin{eqnarray*}
v(t)\leq e^{{\cal I}^\alpha_{t_0}f(t)},~~~t\in[t_0,t^*).
\end{eqnarray*}
Subsequently, from \eqref{ua1sf} we prove
the first inequality of estimate \eqref{uaeIf}.
Because of the continuity of $f$ in $[t_0,t^*)$,
$f$ is bounded in $[t_0,t]$ for each $t\in[t_0,t^*)$.
Simplifying $a(t)e^{{\cal I}^\alpha_{t_0}f(t)}$,
we gain
\begin{eqnarray*}
a(t)e^{{\cal I}^\alpha_{t_0}f(t)}
\!\!\!\!&\leq\!\!\!\!& a(t)\exp(\sup\limits_{s\in[t_0,t]}f(s){\cal I}^\alpha_{t_0}1)
\\
\!\!\!\!&\leq\!\!\!\!& a(t)\exp(\sup\limits_{s\in[t_0,t]}f(s)\frac{|t|^\alpha+|t_0|^\alpha}{\alpha}),
~~~t\in[t_0,t^*),
\end{eqnarray*}
following from Definition \ref{MLfunc} that
the second inequality of estimate \eqref{uaeIf} is true.
Therefore, Lemma \ref{XEE} is proved.
$\qquad\Box$

\subsection{Well-posedness of solutions}

Subsequently, we study the existence and uniqueness,
continuous dependence on initial data
of solutions and continuation of solutions
for the general CFDEs. Consider the initial value problem (IVP) as follows
\begin{eqnarray}
\left\{
\begin{array}{ll}
{\cal T}^\alpha x(t)=f(t,x(t)),~~~(t,x)\in \mathbb{R}^{n+1},
\\
x(t_0)=x_0.
\end{array}
\right.
\label{fcfde-cw}
\end{eqnarray}
Given constants $a,b>0$ and domains
\begin{eqnarray*}
D_+=\{(t,x)\in \mathbb{R}^{n+1}:t\in[t_0-a,t_0+a]\cap\mathbb{R}_+,\|x-x_0\|\leq b\},
~~~t_0\geq0,
\\
D_-=\{(t,x)\in \mathbb{R}^{n+1}:t\in[t_0-a,t_0+a]\cap\mathbb{R}_-,\|x-x_0\|\leq b\},
~~~t_0\leq0,
\end{eqnarray*}
assume that the function $f$ satisfies:
\begin{description}
  \item[(B1)] $f\in C(D_+,\mathbb{R}^n)$ (resp. $f\in C(D_-,\mathbb{R}^n)$);
  \item[(B2)] $f(t,x)$ satisfies Lipschitz condition with respect to $x$ in $D_+$ (resp. $D_-$),
  i.e., there is a positive constant $L$ such that
  \begin{eqnarray*}
    \|f(t,x_1)-f(t,x_2)\|\leq L \|x_1-x_2\|,
    ~~~(t,x_1),(t,x_2)\in D_+~({\rm resp.}~D_-).
  \end{eqnarray*}
\end{description}

\begin{thm}
Suppose that {\bf (B1)} and {\bf (B2)} hold.
Then IVP \eqref{fcfde-cw} has a unique continuous solution in $I_+:=[t_0-\delta_+,t_0+\delta_+]\cap\mathbb{R}_+$
for $t_0\geq0$
(resp. $I_-:=[t_0-\delta_-,t_0+\delta_-]\cap\mathbb{R}_-$ for $t_0\leq0$),
where
\begin{eqnarray*}
\delta_+\!\!\!\!\!\!\!\!&&:=\min\Big{\{}a,\frac{b}{M_+}t_0^{1-\alpha}\Big{\}},
~~~
M_+:=\max_{(t,x)\in D_+}\|f(t,x)\|,
\\
\delta_-\!\!\!\!\!\!\!\!&&:=\min\Big{\{}a,\frac{b}{M_-}(-t_0)^{1-\alpha}\Big{\}},
~~~
M_-:=\max_{(t,x)\in D_-}\|f(t,x)\|.
\end{eqnarray*}
\label{fxxcz-cwa}
\end{thm}

{\bf Proof.}
For convenience, we only discuss the case of $t\in[t_0,t_0+\delta_+]$
for $t_0\geq0$,
all the others can be proved analogously.

Step 1.
By proposition \ref{cfdejfgs},
the equivalent integral equation of IVP \eqref{fcfde-cw} is
\begin{eqnarray}
x(t)=x_0+{\cal I}^\alpha_{t_0}f(t,x(t)),~~~\forall t\in[t_0,t_0+\delta].
\label{jffc}
\end{eqnarray}
Construct the Picard iteration sequence $\{\varphi_n(t)\}$ successively as follows:
\begin{eqnarray}
\left\{
\begin{aligned}
\varphi_0(t)&:=x_0,
\\
\varphi_n(t)&:=x_0+{\cal I}^\alpha_{t_0}f(t,\varphi_{n-1}(t)),~~~n\in\mathbb{N},
~t\in[t_0,t_0+\delta].
\end{aligned}
\right.
\label{pkzbbj}
\end{eqnarray}
We claim that $\varphi_n\in C([t_0,t_0+\delta],\mathbb{R}^n)$ is well-defined
and satisfies
\begin{eqnarray}
\label{vnx0b}
\|\varphi_n(t)-x_0\|\leq b.
\end{eqnarray}
It is obvious that the assertion is true for $n=0$.
Assume that it is also true for $n=k$, then $\varphi_{k+1}\in C([t_0,t_0+\delta],\mathbb{R}^n)$
is also well-defined from the inductive hypothesis
and formula \eqref{pkzbbj}. By {\bf (B1)} we compute
\begin{eqnarray}
\label{vk1x}
\|\varphi_{k+1}(t)-x_0\|\!\!\!\!\!\!\!&&
\leq\Big{\|}{\cal I}^\alpha_{t_0}f(t,\varphi_k(t))\Big{\|}
\nonumber\\
&&\leq M{\cal I}^\alpha_{t_0}1
=\frac{M}{\alpha}(t^\alpha-t_0^\alpha),
~~~t\in[t_0,t_0+\delta].
\end{eqnarray}
By the Lagrange's Mean Value Theorem (\cite{R76}),
there exists a constant $\xi\in(t_0,t)$ such that
\begin{eqnarray*}
|t^\alpha-t_0^\alpha|
=|\alpha \xi^{\alpha-1}(t-t_0)|
\leq \alpha t_0^{\alpha-1}\delta.
\end{eqnarray*}
Substituting the inequality above into \eqref{vk1x},
we realize that \eqref{vnx0b} holds for $n=k+1$.
Then the assertion is true by induction.

Step 2.
We need explain the uniform convergence of sequence $\{\varphi_n(t)\}$ with respect to $t\in[t_0,t_0+\delta]$,
that is the uniform convergence of series
\begin{eqnarray}
\varphi_0(t)+\sum\limits_{k=1}^{+\infty}[\varphi_k(t)-\varphi_{k-1}(t)],~~~\forall t\in[t_0,t_0+\delta],
\label{hsjs}
\end{eqnarray}
since $\varphi_0(t)+\sum\limits_{k=1}^n[\varphi_k(t)-\varphi_{k-1}(t)]=\varphi_n(t)$.
We claim that
\begin{eqnarray}
\label{vkvk1}
\|\varphi_n(t)-\varphi_{n-1}(t)\|\leq\frac{L^{n-1}b^n}{M^{n-1}},
~~~\forall t\in[t_0,t_0+\delta].
\end{eqnarray}
By {\bf (B2)} and induction, one can verify that
\begin{eqnarray*}
\|\varphi_{n+1}(t)-\varphi_n(t)\|\!\!\!\!\!\!\!&&
\leq\Big{\|}{\cal I}^\alpha_{t_0}[f(t,\varphi_n(t))-f(t,\varphi_{n-1}(t))]\Big{\|}\\
&&\leq L{\cal I}^\alpha_{t_0}\|\varphi_n(t)-\varphi_{n-1}(t)\|
\leq L\frac{L^{n-1}b^n}{M^{n-1}}{\cal I}^\alpha_{t_0}1\\
&&\leq \frac{L^{n}b^n}{M^{n-1}}\frac{b}{M}=\frac{L^nb^{n+1}}{M^n},
\end{eqnarray*}
which yields assertion \eqref{vkvk1} is true for $n\in\mathbb{N}$.
By the ratio test, we know that the series $\sum\limits_{k=1}^{+\infty}\cfrac{L^{k-1}b^k}{M^{k-1}}$
is convergent if $b<M/L$, which directly implies that series \eqref{hsjs}
and sequence $\{\varphi_n(t)\}$ are both uniformly convergent in $[t_0,t_0+\delta]$.

Step 3.
Let $\varphi(t):=\lim\limits_{n\to+\infty}\varphi_n(t)$,
then we can easily know that $\varphi\in C([t_0,t_0+\delta],\mathbb{R}^n)$ and fulfills $\|\varphi(t)-x_0\|\leq b$.
By {\bf (B2)} and the continuity of $f$, we obtain $f(t,\varphi(t))=\lim\limits_{n\to+\infty}f(t,\varphi_n(t))$.
Then \eqref{pkzbbj} implies
\begin{eqnarray*}
\lim\limits_{n\to+\infty}\varphi_n(t)=x_0+{\cal I}^\alpha_{t_0}\lim\limits_{n\to+\infty}f(t,\varphi_{n-1}(t)),
\end{eqnarray*}
that is $x(t)=\varphi(t)$ fulfills \eqref{jffc}.
Thus, $\varphi(t)$ is a solution of IVP \eqref{fcfde-cw}.

Step 4.
Suppose that other solution $x(t)=\psi(t)$ also satisfies \eqref{jffc}
associated with initial condition $\psi(t_0)=y_0$,
then by {\bf (B2)} we have
\begin{eqnarray*}
\|\varphi(t)-\psi(t)\|\!\!\!\!\!\!\!&&
\leq\|\varphi(t_0)-\psi(t_0)\|+\Big{\|}{\cal I}^\alpha_{t_0}[f(t,\varphi(t))-f(t,\psi(t))]\Big{\|}\\
&&\leq\|x_0-y_0\|+ L{\cal I}^\alpha_{t_0}\|\varphi(t)-\psi(t)\|,
~~~\forall t\in[t_0,t_0+\delta].
\end{eqnarray*}
We employ Lemma \ref{XEE} to receive
\begin{eqnarray}
\|\varphi(t)-\psi(t)\|\leq \|x_0-y_0\|E_\alpha(L,t)E_\alpha(L,t_0),
~~~\forall t\in[t_0,t_0+\delta].
\label{phipsi}
\end{eqnarray}
Since $E_\alpha(L,t)E_\alpha(L,t_0)$ is bounded for any $t\in[t_0,t_0+\delta]$,
i.e., there exits a constant $M>0$ such that $E_\alpha(L,t)E_\alpha(L,t_0)\leq M$.
For any constant $\varepsilon>0$, if $\|x_0-y_0\|<\varepsilon/M$,
it follows from \eqref{phipsi} that $\|\varphi(t)-\psi(t)\|<\varepsilon$
for $t\in[t_0,t_0+\delta]$.
This proves the continuous dependence on initial data of solutions for IVP \eqref{fcfde-cw}.
Especially, if $x_0=y_0$ then $\varphi(t)\equiv \psi(t)$ for all $t\in[t_0,t_0+\delta]$.
Thus, the uniqueness of solutions is proved.

In conclusion, the proof of Theorem \ref{fxxcz-cwa} is completed.
$\qquad\Box$

\begin{lm}
All solutions of \eqref{TxAx} have maximal interval $\mathbb{R}$.
\end{lm}

{\bf Proof.}
Without loss of generality, we only discuss the case of $t\geq t_0$,
the other case can be proved similarly.
Associated with initial time $t_0$, the equivalent integral equation of \eqref{TxAx} is
\begin{eqnarray*}
x(t)=x(t_0)+{\cal I}^\alpha_{t_0}A(t)x(t),~~~\forall t\in[t_0,+\infty),
\end{eqnarray*}
which implies
\begin{eqnarray*}
\|x(t)\|\leq\|x(t_0)\|+{\cal I}^\alpha_{t_0}\|A(t)\|\|x(t)\|,~~~\forall t\in[t_0,+\infty).
\end{eqnarray*}
It follows from Lemma \ref{XEE} that
\begin{eqnarray*}
\|x(t)\|\leq\|x(t_0)\|E_\alpha\Big{(}\sup\limits_{s\in[t_0,t]}\|A(s)\|,|t|\Big{)}
E_\alpha\Big{(}\sup\limits_{s\in[t_0,t]}\|A(s)\|,|t_0|\Big{)},
~~~\forall t\in[t_0,+\infty),
\end{eqnarray*}
where note that the matrix function $A\in C([t_0,+\infty),\mathbb{R}^{n\times n})$.
Therefore, the solution $x$ exists in interval $[t_0,+\infty)$.
$\qquad\Box$

\subsection{Structure of solutions}

In this subsection, we will present that the linear CFDEs
possess the same structure of solutions as the linear ODEs.
Recall the linear CFDE \eqref{TxAx}.

\begin{pro}
If $x_1,x_2:\mathbb{R}\to \mathbb{R}^n$ are both solutions of \eqref{TxAx},
then $a_1x_1+a_2x_2$ is also a solution of \eqref{TxAx} for any $a_1,a_2\in \mathbb{R}$.
And the set of all solutions of \eqref{TxAx} is an $n$-D linear space.
\label{xxx}
\end{pro}

{\bf Sketch of proof.} One can easily verify the first result of Proposition \ref{xxx}
via property {\bf(A2)}.
Like the corresponding proof in ODEs, see e.g. \cite[p.59]{BV10},
one can verify that solutions $x_1(t),...,x_n(t)$ associated with
initial values $\varepsilon_1,\varepsilon_2,...,\varepsilon_n$
is a basis of $n$-D linear space, if $\varepsilon_1,\varepsilon_2,...,\varepsilon_n$
is a basis of $\mathbb{R}^n$.
$\qquad\Box$

\begin{rmk}
$n\times n$ matrix function $X(t)$, consisting of
$n$ linearly independent solutions $x_1(t),...,x_n(t)$ as its columns,
is also called a {\it fundamental solution} of \eqref{TxAx}.
And for different fundamental solutions $X(t)$ and $Y(t)$,
they can be linearly represented by each other, i.e.,
there exists an invertible linear transformation $C$
such that $Y(t)=X(t)C$ for all $t\in\mathbb{R}$.
\end{rmk}

The proof of this remark is analogous to the case in linear ODEs.

\begin{pro}
The general solution of \eqref{TxAx} associated with initial data $x_0$
can be written as
\begin{eqnarray}
x(t)=X(t)X^{-1}(t_0)x_0,~~~t\in\mathbb{R},
\label{xXXx0}
\end{eqnarray}
where $X(t)$ is any fundamental solution of \eqref{TxAx}.
\label{fzzjbj}
\end{pro}

{\bf Proof.}
Note that the $X(t_0)$ is a nonsingular matrix,
since all its columns, as the initial values of \eqref{TxAx},
form a basis of $\mathbb{R}^n$.
Then \eqref{xXXx0} is well defined, and
$x(t_0)=x_0$. It follows from \eqref{xXXx0} and \eqref{TxAx} that
\begin{eqnarray*}
{\cal T}^\alpha x(t)={\cal T}^\alpha X(t) X^{-1}(t_0)x_0=A(t)X(t)X^{-1}(t_0)x_0=A(t)x(t),
\end{eqnarray*}
that is \eqref{xXXx0} is the general solution of \eqref{TxAx} with $x_0$.
Therefore, Proposition \ref{fzzjbj} is proved.
$\qquad\Box$

\begin{pro}
If $X(t)$ is a fundamental solution of \eqref{TxAx}, then
\begin{eqnarray*}
\det X(t)=\det X(t_0)\exp\Big{(}{\cal I}^\alpha_{t_0}{\rm tr} A(t)\Big{)},
~~~t\in\mathbb{R}.
\end{eqnarray*}
\end{pro}

{\bf Proof.}
Taking $X(t):=[x_{i1}(t),...,x_{in}(t)]$ and $A(t):=[a_{i1}(t),...,a_{in}(t)]$
for $i=1,...,n$, calculate that
\begin{eqnarray*}
{\cal T}^\alpha \det X(t)\!\!\!\!\!\!\!&&=\sum\limits_{i=1}^n\det[{\cal T}^\alpha x_{i1}(t),...,{\cal T}^\alpha x_{in}(t)]\\
&&=\sum\limits_{i=1}^n\det\Big{[}\sum\limits_{j=1}^n a_{ij}(t)x_{j1}(t),\cdots,\sum\limits_{j=1}^n a_{ij}(t)x_{jn}(t)\Big{]}\\
&&=\sum\limits_{i=1}^n a_{ii}(t)\det[x_{i1}(t),\cdots,x_{in}(t)]={\rm tr}A(t)\det X(t).
\end{eqnarray*}
Then $\det X(t)$ is a solution of ${\cal T}^\alpha y(t)={\rm tr}A(t)y(t)$.
From Proposition \ref{TaD}, we derive
\begin{eqnarray*}
y(t)=y(t_0)\exp\Big{(}{\cal I}^\alpha_{t_0}{\rm tr} A(t)\Big{)},
~~~t\in\mathbb{R},
\end{eqnarray*}
yielding the desired result.
$\qquad\Box$

\subsection{Autonomous systems}

As a special case of \eqref{TxAx}, here we consider the linear autonomous system
\begin{eqnarray}
{\cal T}^\alpha x=Ax,~~~(t,x)\in\mathbb{R}^{n+1},
\label{TxtAx}
\end{eqnarray}
where $A$ is an $n\times n$ real constant matrix.
To solving \eqref{TxtAx}, we first provide the concept of the Mittag-Leffler-type of a matrix.

\begin{df}
The Mittag-Leffler-type of an $n\times n$ real constant matrix $A$ is defined as
\begin{eqnarray}
E_\alpha(A,1):=\sum_{k=0}^{+\infty}\frac{A^k}{\alpha^kk!},
\label{ejs}
\end{eqnarray}
and denote $E_\alpha(0,1)=I$ for convention.
\end{df}

\begin{pro}
The power series in \eqref{ejs} is convergent for any matrix $A$.
\end{pro}

{\bf Proof.}
In deed, it is obvious that $\|A^k\|\leq \|A\|^k$ for all $k\in \mathbb{N}$, so
\begin{eqnarray*}
\sum_{k=0}^{+\infty}\frac{\|A^k\|}{\alpha^kk!}\leq \sum_{k=0}^{+\infty}\frac{\|A\|^k}{\alpha^kk!}=E_\alpha(\|A\|,1)<+\infty.
\end{eqnarray*}
It yields that the series in \eqref{ejs} is convergent.
$\qquad\Box$

Recall the Jordan canonical form in ODEs as follows
\begin{eqnarray}
\label{P1AP}
P^{-1}AP
:=
\left[ {\begin{array}{ccc}
J_1 & \cdots& 0 \\
0&\ddots &0 \\
0 &\cdots &J_l \\
\end{array}} \right],
~~~
J_i=
\left[ {\begin{array}{ccc}
\lambda_i & 1 & 0 \\
 0& \ddots & 1 \\
 0 & 0& \lambda_i\\
\end{array}} \right],~~~i=1,2,...,l,
\end{eqnarray}
where $P$ is an $n\times n$ nonsingular complex matrix
and $\lambda_i$ is an eigenvalue of $A$.

Thus, the Mittag-Leffler-type of a matrix can be easily computed as follows.
\begin{pro}
Let $A$ is an $n\times n$ real matrix with Jordan canonical form in \eqref{P1AP},
then
\begin{eqnarray*}
E_\alpha(A,1)=PE_\alpha(P^{-1}AP,1)P^{-1}=
P\left[ {\begin{array}{ccc}
E_\alpha(J_1,1) & \cdots& 0 \\
0&\ddots &0 \\
0 &\cdots &E_\alpha(J_l,1) \\
\end{array}} \right]P^{-1}.
\end{eqnarray*}
\label{MLTm}
\end{pro}

{\bf Proof.}
Since $(P^{-1}AP)^k=P^{-1}A^kP$ for any $k\in \mathbb{N}$, one can compute
\begin{eqnarray*}
E_\alpha(P^{-1}AP,1)=P^{-1}\Big{(}\sum\limits_{k=0}^{+\infty}\frac{A^k}{\alpha^k k!}\Big{)}P=P^{-1}E_\alpha(A,1)P.
\end{eqnarray*}
Subsequently, one can easily verify that
\begin{eqnarray*}
E_\alpha(P^{-1}AP,1)=
\sum\limits_{k=0}^{+\infty}\frac{1}{\alpha^k k!}
\left[ {\begin{array}{ccc}
J_1^k & \cdots& 0 \\
0&\ddots &0 \\
0 &\cdots &J_l^k \\
\end{array}} \right]=
\left[ {\begin{array}{ccc}
E_\alpha(J_1,1) & \cdots& 0 \\
0&\ddots &0 \\
0 &\cdots &E_\alpha(J_l,1) \\
\end{array}} \right],
\end{eqnarray*}
which implies the required conclusion.
$\qquad\Box$

Further, one can verify the following formula.

\begin{pro}
If $A=\lambda I+N$, where the nilpotent matrix $N$ is
\begin{eqnarray*}
N=
\left[ {\begin{array}{ccc}
0 & 1& 0 \\
0&\ddots &1 \\
0 &0 &0 \\
\end{array}} \right],
\end{eqnarray*}
then the following expression fulfills:
\begin{eqnarray*}
E_\alpha(A,1)=E_\alpha(\lambda,1)\Big{(}I+\frac{N}{\alpha}+\frac{N^2}{\alpha^2 2!}+\cdots+\frac{N^{n-1}}{\alpha^{(n-1)} (n-1)!}\Big{)}.
\end{eqnarray*}
\end{pro}

Relying on the preliminaries above, one can solve \eqref{TxtAx} as follows.

\begin{lm}
The matrix $E_\alpha(A,t)$ is a fundamental solution of \eqref{TxtAx}
for all $t\in\mathbb{R}$.
\label{TMXa}
\end{lm}

{\bf Proof.}
Note that
\begin{eqnarray*}
E_\alpha(A,t):=
\left\{
\begin{array}{ll}
\sum\limits_{k=0}^{+\infty}\cfrac{A^k t^{\alpha k}}{\alpha ^kk!},
~~~t\in\mathbb{R}_+,
\\
\sum\limits_{k=0}^{+\infty}\cfrac{(-A)^k(-t)^{\alpha k}}{\alpha ^kk!},
~~~t\in\mathbb{R}_-.
\end{array}
\right.
\end{eqnarray*}
The following calculus implies that power series $E_\alpha(A,t)$ is convergent for all $t\in\mathbb{R}$:
\begin{eqnarray*}
\sum_{k=0}^{+\infty}\frac{|t|^{\alpha k}\|A^k\|}{\alpha^{k}k!}
\leq \sum_{k=0}^{+\infty}\frac{|t|^{\alpha k} \|A\|^k}{\alpha^kk!}=E_\alpha(\|A\|,|t|)<+\infty.
\end{eqnarray*}
And one can compute that
\begin{eqnarray*}
{\cal T}^\alpha E_\alpha(A,t)\!\!\!\!\!\!&&=\sum_{k=1}^{+\infty}\frac{t^{1-\alpha} t^{\alpha k-1}A^k}{\alpha^{k-1}(k-1)!}
\\
\!\!\!\!\!\!&&=A\sum_{k=1}^{+\infty}\frac{t^{\alpha(k-1)}A^{k-1}}{\alpha^{k-1}(k-1)!}=AE_\alpha(A,t),
~~~t\in\mathbb{R}_+,
\\
{\cal T}^\alpha E_\alpha(A,t)\!\!\!\!\!\!&&=\sum_{k=1}^{+\infty}\frac{-(-t)^{1-\alpha} (-t)^{\alpha k-1}(-A)^k}{\alpha^{k-1}(k-1)!}
\\
\!\!\!\!\!\!&&=A\sum_{k=1}^{+\infty}\frac{(-t)^{\alpha(k-1)}(-A)^{k-1}}{\alpha^{k-1}(k-1)!}=AE_\alpha(A,t),
~~~t\in\mathbb{R}_-,
\end{eqnarray*}
which yields the wanted result.
$\qquad\Box$

Both Proposition \ref{fzzjbj} and Lemma \ref{TMXa} lead to
the following result on general solutions of \eqref{TxtAx}.

\begin{pro}
The general solution of \eqref{TxtAx} associated with initial data $x_0$
can be expressed as
\begin{eqnarray*}
x(t)=\frac{E_\alpha(A,t)}{E_\alpha(A,t_0)}x_0,
~~~t\in\mathbb{R}.
\end{eqnarray*}
\label{XEAEA}
\end{pro}

\subsection{Variation of constants formula}

In this subsection, we focus on the perturbation of linear CFDE \eqref{TxAx}.
Consider the inhomogeneous linear CFDE
\begin{eqnarray}
{\cal T}^\alpha x=A(t)x+f(t),~~~(t,x)\in\mathbb{R}^{n+1},
\label{TxAxf}
\end{eqnarray}
where $f\in C(\mathbb{R},\mathbb{R}^{n})$ and matrix function
$A\in C(\mathbb{R},\mathbb{R}^{n\times n})$.

One can verify the analogous properties to ODEs as follows.

\begin{pro}
Like ODEs, if both $x^*_1(t)$ and $x^*_2(t)$ are solutions of \eqref{TxAxf},
then $x^*_1(t)-x^*_2(t)$ is a solution of \eqref{TxAx}.
On the other hand, if $x(t)$ and $x^*(t)$ are solutions of \eqref{TxAx} and \eqref{TxAxf} respectively,
then $x(t)+x^*(t)$ is also a solution of \eqref{TxAxf}.
\end{pro}

These properties can easily lead to the following structure of general solutions for \eqref{TxAxf}.

\begin{pro}
If $x^*(t)$ is a solution of \eqref{TxAxf},
then general solutions of \eqref{TxAxf} associated with initial data $x_0$
can be represented as
\begin{eqnarray*}
x(t)=X(t)X^{-1}(t_0)x_0+x^*(t),~~~t\in\mathbb{R},
\end{eqnarray*}
where $X(t)$ is any fundamental solution of \eqref{TxAx}.
\end{pro}

Next, we will give the variation of constants formula for \eqref{TxAxf}.

\begin{thm}
Let $X(t)$ is a fundamental matrix of \eqref{TxAx},
then the general solutions of \eqref{TxAxf} associated with initial data $x_0$
can be given by
\begin{eqnarray}
x(t)=X(t)X^{-1}(t_0)x_0+X(t){\cal I}^\alpha_{t_0}X^{-1}(t)f(t),
~~~t\in\mathbb{R}.
\label{csbygs}
\end{eqnarray}
Particularly, if $A(t)$ degenerates into an $n\times n$ real constant matrix $A$,
the variation of constants formula \eqref{csbygs} becomes the form
\begin{eqnarray*}
x(t)=\frac{E_\alpha(A,t)}{E_\alpha(A,t_0)}x_0+E_\alpha(A,t){\cal I}^\alpha_{t_0}\frac{f(t)}{E_\alpha(A,t)},
~~~t\in\mathbb{R}.
\end{eqnarray*}
\label{fzzcsbygs}
\end{thm}

{\bf Proof.}
By Proposition \ref{fzzjbj}, the general solutions of \eqref{TxAx}
is $x(t)=X(t)X^{-1}(t_0)c$ for any constant column vector $c$.
By analogy with ODEs, let
\begin{eqnarray}
x(t)=X(t)X^{-1}(t_0)c(t)
\label{xXX-1c}
\end{eqnarray}
is the general solutions of \eqref{TxAxf},
then
\begin{eqnarray}
{\cal T}^\alpha x(t)\!\!\!\!\!\!\!\!&&=({\cal T}^\alpha X(t))X^{-1}(t_0)c(t)+X(t)X^{-1}(t_0){\cal T}^\alpha c(t)
\nonumber \\
&&=A(t)X(t)X^{-1}(t_0)c(t)+f(t).
\label{TxTXX}
\end{eqnarray}
Since ${\cal T}^\alpha X(t)=A(t)X(t)$, it follows from \eqref{TxTXX} that
\begin{eqnarray}
X(t)X^{-1}(t_0){\cal T}^\alpha c(t)=f(t),
~~~t\in\mathbb{R}.
\label{XX-1T}
\end{eqnarray}
Integrating both hands sides of \eqref{XX-1T} from $t_0$ to $t$,
we receive
\begin{eqnarray}
c(t)=c(t_0)+{\cal I}^\alpha_{t_0}X(t_0)X^{-1}(t)f(t),
~~~t\in\mathbb{R}.
\label{ccIXX}
\end{eqnarray}
Substituting \eqref{ccIXX} into \eqref{xXX-1c},
we receive
\begin{eqnarray*}
x(t)\!\!\!\!\!\!\!&&=X(t)X^{-1}(t_0)\Big{[}c(t_0)+{\cal I}^\alpha_{t_0}X(t_0)X^{-1}(t)f(t)\Big{]}\\
&&=X(t)X^{-1}(t_0)c(t_0)+X(t){\cal I}^\alpha_{t_0}X^{-1}(t)f(t),
~~~t\in\mathbb{R},
\end{eqnarray*}
where \eqref{xXX-1c} guarantees that $c(t_0)=x_0$.
It yields that \eqref{csbygs} holds.
The special case $X(t)=E_\alpha(A,t)$ can be obtained naturally.
Thus, Theorem \ref{fzzcsbygs} is proved.
$\qquad\Box$

\begin{pro}
If an $n\times n$ real constant matrix $A$ has only eigenvalues
with negative real part,
then there exist constants $K,\lambda>0$ such that
\begin{eqnarray*}
\|E_\alpha(A,t)\|\leq KE_\alpha(-\lambda,t),
~~~t\in\mathbb{R}_+.
\end{eqnarray*}
\end{pro}

The proof is similar to the case in ODEs, referred Proposition 2.27 in \cite[p.77]{BV10}.

\section{Stability and Mittag-Leffler dichotomy}
\setcounter{equation}{0}
In this section, we study the concepts of stability and Mittag-Leffler dichotomy of CFDEs.
Before this, Souahi, Makhlouf and Hammami(\cite{SMH17}) combined Lyapunov stability and properties of conformable fractional
derivative given by Abdeljawad(\cite{A15})
to raise the concepts of stability, asymptotic stability and
fractional exponential stability for the nonlinear system \eqref{fcfde-cw}.
For the nonautonomous linear CFDE \eqref{TxAx},
the definitions of uniform stability and uniformly asymptotic stability are more essential.

Based on the definition of stability for CFDEs described in \cite{SMH17},
we introduce the following definition of uniformly stability
analogous to the corresponding concept of ODEs in e.g. \cite[p.1]{C78}.

\begin{df}
The solution $\hat{x}(t)$ of system \eqref{fcfde-cw} is said to be
\begin{description}
  \item[(C1)] uniformly stable, if for any $\varepsilon>0$ there exists $\delta:=\delta(\varepsilon)>0$
  such that any solution $x(t)$ of \eqref{fcfde-cw} satisfies for some $s\geq0$,
  the inequality $\|x(s)-\hat{x}(s)\|<\delta$ implies $\|x(t)-\hat{x}(t)\|<\varepsilon$ for all $t\geq s$;
  \item[(C2)] attractive, if there exists $\delta_0>0$ and $T:=T(\varepsilon)>0$ for any $\varepsilon>0$
          such that for some $s\geq0$, the inequality $\|x(s)-\hat{x}(s)\|<\delta_0$ implies $\|x(t)-\hat{x}(t)\|<\varepsilon$ for all $t\geq s+T$;
  \item[(C3)] uniformly asymptotically stable, if it is uniformly stable and attractive.
\end{description}
\label{yzwdx}
\end{df}

The following definition is on the {\it Mittag-Leffler stability}.

\begin{df}
The solution $x_*=0$ of system \eqref{fcfde-cw} is Mittag-Leffler stable if
\begin{eqnarray*}
\|x(t)\|\leq K\frac{E_\alpha(\lambda,t_0)}{E_\alpha(\lambda,t)}\|x_0\|,~~~t\geq t_0,
\end{eqnarray*}
where constants $K,\lambda>0$.
\label{xKE}
\end{df}

More generally, we focus on the significant application of Definition \ref{yzwdx} to linear equation \eqref{TxAx}.

\begin{pro}
Suppose $X(t)$ is a fundamental matrix of \eqref{TxAx} and $c$ is a real constant,
then solution $x_*=0$ of \eqref{TxAx} is said to be
\begin{description}
  \item[(D1)] stable for any $t_0\in \mathbb{R}$ if and only if there exists $K:=K(t_0)>0$ such that
  \begin{eqnarray*}
   \|X(t)\|\leq K,~~~t_0\leq t<+\infty;
   \end{eqnarray*}
  \item[(D2)] uniformly stable for $t_0\geq c$ if and only if there exists $K:=K(c)>0$ such that
   \begin{eqnarray*}
   \|X(t)X^{-1}(s)\|\leq K,~~~t_0\leq s\leq t <+\infty;
   \end{eqnarray*}
  \item[(D3)] asymptotically stable for any $t_0\in \mathbb{R}$ if and only if $\lim\limits_{t\rightarrow+\infty}\|X(t)\|=0$;
  \item[(D4)] uniformly asymptotically stable for $t_0\geq c$ if and only if there exist $K:=K(c)>0$
  and $\lambda:=\lambda(c)>0$ such that
  \begin{eqnarray}
   \|X(t)X^{-1}(s)\|\leq K\frac{E_\alpha(\lambda,s)}{E_\alpha(\lambda,t)},~~~t_0\leq s\leq t<+\infty.
   \label{XX-1KE}
   \end{eqnarray}
\end{description}
Particularly, {\bf (D1)-(D4)} all hold for autonomous system \eqref{TxtAx},
if fundamental matrix $X(t)$ is replaced by $E_\alpha(A,t)$.
\label{wdx}
\end{pro}

The proof of {\bf (D1)}-{\bf (D4)} can refer to the Theorem 2.1 in \cite[p.84]{H69}.
In particular, since the Mittag-Leffler stability implies the uniformly asymptotic stability,
one can simply verify conclusion {\bf (D4)}.
The definitions of the corresponding stabilities above for Caputo FDEs
had been proposed in references e.g. \cite[p.140]{DMD19}.

Next, we shall propose the concept of {\it Mittag-Leffler dichotomy} for linear CFDE \eqref{TxAx}.

\begin{df}
Suppose that $X(t)$ is a fundamental matrix of \eqref{TxAx}.
The equation \eqref{TxAx} possesses a Mittag-Leffler dichotomy
if there exists a projection matrix $P$, i.e. $P^2=P$, and positive constants $N_i$, $\beta_i$ $(i=1,2)$ such that
\begin{eqnarray}
\begin{split}
\|X(t)PX^{-1}(s)\|\leq N_1\frac{E_\alpha(\beta_1,s)}{E_\alpha(\beta_1,t)},~~~t\geq s,\\
\|X(t)(I-P)X^{-1}(s)\|\leq N_2\frac{E_\alpha(\beta_2,t)}{E_\alpha(\beta_2,s)},~~~s\geq t.
\end{split}
\label{XPXK}
\end{eqnarray}
In particular, \eqref{TxAx} possesses an ordinary dichotomy if \eqref{XPXK} hold with $\beta_1=\beta_2=0$.
\label{fzzer}
\end{df}

Finally, we concern perturbation of nonautonomous linear CFDE \eqref{TxAx}.
Consider the perturbed equation
\begin{eqnarray}
{\cal T}^\alpha x=A(t)x+f(t,x),~~~(t,x)\in\mathbb{R}^{n+1},
\label{TxAxft}
\end{eqnarray}
where $f\in C(\mathbb{R}^{n+1},\mathbb{R}^{n})$ and matrix function $A\in C(\mathbb{R},\mathbb{R}^{n\times n})$.

The following conclusion give out the projection form of equivalent integral equation
and the existence of bounded solutions for equation \eqref{TxAxft}.

\begin{lm}
Suppose that function $f\in C(\mathbb{R}^{n+1},\mathbb{R}^{n})$,
$P$ is a projection matrix given in Definition \ref{fzzer} and
equation \eqref{TxAx} possesses a Mittag-Leffler dichotomy.
If $x\in C_b([t_0,+\infty),\mathbb{R}^n)$ is a solution of \eqref{TxAxft} with $x(t_0)=x_0$
for constant $t_0\in \mathbb{R}_+$, then
\begin{eqnarray}
\begin{split}
x(t)=&X(t)PX^{-1}(t_0)x_0+X(t)P{\cal I}^\alpha_{t_0}X^{-1}(t)f(t,x(t))\\
&+X(t)(I-P){\cal I}^\alpha_{+\infty}X^{-1}(t)f(t,x(t)),
~~~t\geq t_0.
\label{xXPX1}
\end{split}
\end{eqnarray}
If $x\in C_b((-\infty, t_0],\mathbb{R}^n)$ is a solution of \eqref{TxAxft} with $x(t_0)=x_0$
for constant $t_0\in \mathbb{R}_-$, then
\begin{eqnarray}
\begin{split}
x(t)=&X(t)(I-P)X^{-1}(t_0)x_0+X(t)(I-P){\cal I}^\alpha_{t_0}X^{-1}(t)f(t,x(t))\\
&+X(t)P{\cal I}^\alpha_{-\infty}X^{-1}(t)f(t,x(t)),
~~~t\leq t_0.
\label{xXPX2}
\end{split}
\end{eqnarray}
Conversely, any bounded solution of \eqref{xXPX1} or \eqref{xXPX2} is a solution of \eqref{TxAxft}.
\label{efds}
\end{lm}

{\bf Proof.}
For convenience, we only prove \eqref{xXPX1},
formula \eqref{xXPX2} can be proved in an analogous manner.
Assume $x(t)$ is a bounded solution of \eqref{TxAxft}
and $M:=\sup\limits_{t\in[t_0,+\infty)}\|x(t)\|$ for $t\geq t_0$.
The continuity of $f$ implies that
there exists a positive constant $N$ such that $N:=\sup\limits_{t\in[t_0,+\infty)}\|f(t,x(t))\|$.
By the variation of constants formula \eqref{csbygs}, for any $\tau\geq t_0$,
the solution $x(t)$ satisfies
\begin{eqnarray}
\begin{split}
X(t)(I-P)X^{-1}(t)x(t)=&X(t)(I-P)X^{-1}(\tau)x(\tau)\\
&+X(t)(I-P){\cal I}^\alpha_{\tau}X^{-1}(t)f(t,x(t)),
~~~t,\tau\geq t_0,
\label{XPX1}
\end{split}
\end{eqnarray}
where the following estimate can be obtained by \eqref{XPXK}
\begin{eqnarray*}
\|X(t)(I-P)X^{-1}(\tau)x(\tau)\|\!\!\!\!\!\!\!\!&&
\leq N_2\frac{E_\alpha(\beta_2,t)}{E_\alpha(\beta_2,\tau)}\sup\limits_{t\in[t_0,+\infty)}\|x(t)\|
\\
\!\!\!\!\!\!\!\!&&\leq MN_2\frac{E_\alpha(\beta_2,t)}{E_\alpha(\beta_2,\tau)},
~~~t,\tau\geq t_0.
\end{eqnarray*}
It yields that
\begin{eqnarray*}
\lim\limits_{\tau\to+\infty}\|X(t)(I-P)X^{-1}(\tau)x(\tau)\|=0,
~~~t\geq t_0.
\end{eqnarray*}
On the other hand, in integral equation \eqref{XPX1} for $t\geq t_0$,
\begin{eqnarray*}
\|X(t)(I-P){\cal I}^\alpha_{\tau}X^{-1}(t)f(t,x(t))\|
\!\!\!\!\!\!\!\!&&\leq
NN_2E_\alpha(\beta_2,t)|{\cal I}^\alpha_{\tau}E_\alpha(-\beta_2,t)|
\\
\!\!\!\!\!\!\!\!&&\leq
NN_2E_\alpha(\beta_2,t)\int^{\tau}_{t}s^{\alpha-1}\exp\Big{(}-\beta_2 \cfrac{s^\alpha}{\alpha}\Big{)}ds
\\
\!\!\!\!\!\!\!\!&&\leq
\frac{NN_2E_\alpha(\beta_2,t)}{\beta_2}(E_\alpha(-\beta_2,t)-E_\alpha(-\beta_2,\tau))
\\
\!\!\!\!\!\!\!\!&&\leq
\frac{NN_2}{\beta_2},
\end{eqnarray*}
which implies that
\begin{eqnarray*}
\|X(t)(I-P){\cal I}^\alpha_{+\infty}X^{-1}(t)f(t,x(t))\|<+\infty,
~~~t\geq t_0.
\end{eqnarray*}
It follows from \eqref{XPX1} that
\begin{eqnarray}
\begin{split}
X(t)(I-P)X^{-1}(t)x(t)=X(t)(I-P){\cal I}^\alpha_{+\infty}X^{-1}(t)f(t,x(t)),
~~~t\geq t_0.
\end{split}
\label{I-PxXI-P}
\end{eqnarray}
From the variation of constants formula \eqref{csbygs},
it also follows that for $t\geq t_0$,
\begin{eqnarray}
\begin{split}
X(t)PX^{-1}(t)x(t)=X(t)PX^{-1}(t_0)x(t_0)+X(t)P{\cal I}^\alpha_{t_0}X^{-1}(t)f(t,x(t)).
\end{split}
\label{PxXPX-1}
\end{eqnarray}
Since $x(t)=X(t)PX^{-1}(t)x(t)+X(t)(I-P)X^{-1}(t)x(t)$,
substituting \eqref{I-PxXI-P} and \eqref{PxXPX-1} into it,
we attain \eqref{xXPX1}.
And the converse conclusion can be verified by direct calculation
to end the proof.
$\qquad\Box$

The following Lemma is the fractional-order version of projected integral inequality.

\begin{lm}
Suppose that $N_i$, $\beta_i$ and $\varepsilon$ are all positive constants for $i=1,2$,
and bounded continuous nonnegative solutions $u(t)$ satisfy
\begin{eqnarray}
\begin{split}
u(t)\leq& N_1\frac{E_\alpha(\beta_1,t_0)}{E_\alpha(\beta_1,t)}
+\frac{\varepsilon N_1}{E_\alpha(\beta_1,t)}{\cal I}^\alpha_{t_0}E_\alpha(\beta_1,t)u(t)
\\
&-\varepsilon N_2E_\alpha(\beta_2,t){\cal I}^\alpha_{+\infty}\frac{u(t)}{E_\alpha(\beta_2,t)},
~~~t\geq t_0\geq0,
\label{xEEN1}
\end{split}
\\
\begin{split}
u(t)\leq& N_2\frac{E_\alpha(\beta_2,t)}{E_\alpha(\beta_2,t_0)}
-\varepsilon N_2 E_\alpha(\beta_2,t){\cal I}^\alpha_{t_0}\frac{u(t)}{E_\alpha(\beta_2,t)}
\\
&+\frac{\varepsilon N_1}{E_\alpha(\beta_1,t)}{\cal I}^\alpha_{-\infty}E_\alpha(\beta_1,t)u(t),
~~~t\leq t_0\leq0.
\label{xEEN2}
\end{split}
\end{eqnarray}
Set that
\begin{eqnarray*}
\theta:=\varepsilon\Big{(}\frac{N_1}{\beta_1}+\frac{N_2}{\beta_2}\Big{)},
~~~
K_i:=\frac{N_i}{1-\theta},~~~\lambda_i:=\beta_i-\frac{\varepsilon N_i}{1-\theta},
~~~
i=1,2.
\nonumber
\end{eqnarray*}
If $\theta<1$, then
\begin{eqnarray*}
u(t)\leq
\left\{
\begin{array}{ll}
K_1\cfrac{E_\alpha(\lambda_1,t_0)}{E_\alpha(\lambda_1,t)},~~~t\geq t_0,
\\
K_2\cfrac{E_\alpha(\lambda_2,t)}{E_\alpha(\lambda_2,t_0)},~~~t\leq t_0.
\end{array}
\right.
\end{eqnarray*}
\label{efbds1}
\end{lm}

{\bf Proof.}
Without loss of generality, we only consider inequality \eqref{xEEN1},
because inequality \eqref{xEEN2} can be changed into \eqref{xEEN1}
through transformations $t\to-t$ and $t_0\to-t_0$.
Next, we need to verify $\lim\limits_{t\rightarrow{+\infty}}u(t)=0$.
In deed, since $u(t)$ is bounded, let $\sigma:=\limsup\limits_{t\rightarrow{+\infty}}u(t)$.
If $\sigma>0$ and for any constant $\vartheta$ satisfying $\theta<\vartheta<1$,
there exists $t_1\geq t_0$ such that for any $t\geq t_1$ we have $u(t)\leq\vartheta^{-1}\sigma$.
For $t\geq t_1$ we compute
\begin{eqnarray*}
u(t)\leq \!\!\!\!\!\!\!\!\!&& N_1\frac{E_\alpha(\beta_1,t_0)}{E_\alpha(\beta_1,t)}
+\frac{\varepsilon N_1}{E_\alpha(\beta_1,t)}{\cal I}^\alpha_{t_0}E_\alpha(\beta_1,t_1)u(t_1)
\\
&&+\vartheta^{-1}\sigma\Big{[}\frac{\varepsilon N_1}{E_\alpha(\beta_1,t)}{\cal I}^\alpha_{t_1}E_\alpha(\beta_1,t)
-\varepsilon N_2E_\alpha(\beta_2,t){\cal I}^\alpha_{+\infty}\frac{1}{E_\alpha(\beta_2,t)}\Big{]}\\
\leq \!\!\!\!\!\!\!\!\!&& N_1\frac{E_\alpha(\beta_1,t_0)}{E_\alpha(\beta_1,t)}
+\frac{\varepsilon N_1}{E_\alpha(\beta_1,t)}{\cal I}^\alpha_{t_0}E_\alpha(\beta_1,t_1)u(t_1)
+\vartheta^{-1}\sigma\varepsilon\Big{(}\frac{ N_1}{\beta_1}+\frac{ N_2}{\beta_2}\Big{)}.
\end{eqnarray*}
Since $\theta<\vartheta<1$, the upper limit of the right hand side of the inequality above is less than $\sigma$ as $t\to +\infty$.
It follows from the inequality above that
\begin{eqnarray*}
\sigma\leq\vartheta^{-1}\sigma\varepsilon\Big{(}\frac{ N_1}{\beta_1}+\frac{ N_2}{\beta_2}\Big{)}<\sigma,
\end{eqnarray*}
that is a contradiction.
Hence, $\sigma=0$ and $\lim\limits_{t\rightarrow{+\infty}}u(t)=0$.

Set
$
v(t):=\sup\limits_{\tau\geq t}u(\tau).
$
Obviously, the function $v(t)$ is nonincreasing
and for any $t\geq t_0$, there exists $t_2\geq t$ such that for $t\leq s\leq t_2$,
$v(t)=u(t_2)=v(s)$.
Replacing $t$ in \eqref{xEEN1} with $t_2$, for $t\geq t_0$ we calculate that
\begin{eqnarray*}
v(t)=u(t_2)\leq\!\!\!\!\!\!\!\!&& N_1\frac{E_\alpha(\beta_1,t_0)}{E_\alpha(\beta_1,t_2)}
+\frac{\varepsilon N_1}{E_\alpha(\beta_1,t_2)}{\cal I}^\alpha_{t_0}E_\alpha(\beta_1,t_2)u(t_2)\\
&&-\varepsilon N_2E_\alpha(\beta_2,t_2){\cal I}^\alpha_{+\infty}\frac{u(t_2)}{E_\alpha(\beta_2,t_2)}
\\
\leq\!\!\!\!\!\!\!\!&& N_1\frac{E_\alpha(\beta_1,t_0)}{E_\alpha(\beta_1,t)}
+\frac{\varepsilon N_1}{E_\alpha(\beta_1,t)}{\cal I}^\alpha_{t_0}E_\alpha(\beta_1,t_2)u(t_2)
\\
&&-\varepsilon N_2E_\alpha(\beta_2,t_2){\cal I}^\alpha_{+\infty}\frac{u(t_2)}{E_\alpha(\beta_2,t_2)}\\
\leq\!\!\!\!\!\!\!\!&& N_1\frac{E_\alpha(\beta_1,t_0)}{E_\alpha(\beta_1,t)}
+\frac{\varepsilon N_1}{E_\alpha(\beta_1,t)}{\cal I}^\alpha_{t_0}E_\alpha(\beta_1,t)u(t)\\
&&-v(t)\Big{[}\frac{\varepsilon N_1}{E_\alpha(\beta_1,t_2)}{\cal I}^\alpha_{t_2}E_\alpha(\beta_1,t)
+\varepsilon N_2E_\alpha(\beta_2,t_2){\cal I}^\alpha_{+\infty}\frac{1}{E_\alpha(\beta_2,t_2)}\Big{]}\\
\leq\!\!\!\!\!\!\!\!&& N_1\frac{E_\alpha(\beta_1,t_0)}{E_\alpha(\beta_1,t)}
+\frac{\varepsilon N_1}{E_\alpha(\beta_1,t)}{\cal I}^\alpha_{t_0}E_\alpha(\beta_1,t)u(t)
+v(t)\varepsilon\Big{(}\frac{ N_1}{\beta_1}+\frac{ N_2}{\beta_2}\Big{)}.
\end{eqnarray*}
Put $w(t):=\cfrac{E_\alpha(\beta_1,t)}{E_\alpha(\beta_1,t_0)}v(t)$, then it follows from the definition of $v$ that
\begin{eqnarray*}
w(t)\!\!\!\!\!\!\!\!&&\leq N_1
+\frac{\varepsilon N_1}{E_\alpha(\beta_1,t_0)}{\cal I}^\alpha_{t_0}E_\alpha(\beta_1,t)v(t)+\theta w(t)\\
&&= N_1+\varepsilon N_1{\cal I}^\alpha_{t_0}w(t)+\theta w(t),
~~~t\geq t_0,
\end{eqnarray*}
that is
\begin{eqnarray*}
w(t)\leq \frac{N_1}{1-\theta}+\frac{\varepsilon N_1}{1-\theta}{\cal I}^\alpha_{t_0}w(t),
~~~t\geq t_0.
\end{eqnarray*}
Applying Lemma \ref{XEE} to the inequality above, we attain
\begin{eqnarray*}
w(t)\leq \frac{N_1}{1-\theta}\exp({\cal I}^\alpha_{t_0}\frac{\varepsilon N_1}{1-\theta})
=\frac{N_1}{1-\theta}\frac{E_\alpha(\cfrac{\varepsilon N_1}{1-\theta},t)}
{E_\alpha(\cfrac{\varepsilon N_1}{1-\theta},t_0)},
~~~t\geq t_0.
\end{eqnarray*}
Combining with the definitions of $v$ and $w$, we acquire
\begin{eqnarray*}
u(t)\!\!\!\!\!\!\!\!&&\leq\frac{N_1}{1-\theta}\frac{E_\alpha(\cfrac{\varepsilon N_1}{1-\theta},t)}
{E_\alpha(\cfrac{\varepsilon N_1}{1-\theta},t_0)}
\frac{E_\alpha(\beta_1,t_0)}{E_\alpha(\beta_1,t)}
=K_1\frac{E_\alpha(\lambda_1,t_0)}{E_\alpha(\lambda_1,t)},~~~t\geq t_0,
\end{eqnarray*}
where $K_1=\cfrac{N_1}{1-\theta}$ and $\lambda_1=\beta_1-\cfrac{\delta N_1}{1-\theta}$.
Therefore, Lemma \ref{efbds1} is proved.
$\qquad\Box$

As a corollary of Lemma \ref{efbds1}, we introduce a more useful result in estimate of dichotomy.

\begin{co}
Suppose that $N_i$, $\beta_i$ and $\varepsilon$ are all positive constants for $i=1,2$,
and bounded continuous nonnegative solutions $u(t)$ satisfy
\begin{eqnarray}
\begin{split}
u(t)\leq& N_2\frac{E_\alpha(\beta_2,t)}{E_\alpha(\beta_2,s)}
+\frac{\varepsilon N_1}{E_\alpha(\beta_1,t)}{\cal I}^\alpha_{t_0}E_\alpha(\beta_1,t)u(t)
\\
&-\varepsilon N_2 E_\alpha(\beta_2,t){\cal I}^\alpha_{s}\frac{u(t)}{E_\alpha(\beta_2,t)},
~~~s\geq t\geq t_0\geq0,
\label{xEEN3}
\end{split}
\\
\begin{split}
u(t)\leq& N_1\frac{E_\alpha(\beta_1,s)}{E_\alpha(\beta_1,t)}
-\varepsilon N_2 E_\alpha(\beta_2,t){\cal I}^\alpha_{t_0}\frac{u(t)}{E_\alpha(\beta_2,t)}
\\
&+\frac{\varepsilon N_1}{E_\alpha(\beta_1,t)}{\cal I}^\alpha_{s}E_\alpha(\beta_1,t)u(t),
~~~s\leq t\leq t_0\leq0.
\label{xEEN4}
\end{split}
\end{eqnarray}
If $\theta<1$, then
\begin{eqnarray*}
u(t)\leq
\left\{
\begin{array}{ll}
K_2\cfrac{E_\alpha(\lambda_2,t)}{E_\alpha(\lambda_2,s)},~~~s\geq t\geq t_0,
\\
K_1\cfrac{E_\alpha(\lambda_1,s)}{E_\alpha(\lambda_1,t)},~~~s\leq t\leq t_0,
\end{array}
\right.
\end{eqnarray*}
where $\theta$, $K_i$ and $\lambda_i$ were all defined in Lemma \ref{efbds1}.
\label{efbds2}
\end{co}

{\bf Proof.}
Without loss of generality, we only consider inequality \eqref{xEEN3},
because inequality \eqref{xEEN4} can be changed into \eqref{xEEN3}
through transformations $t\to-t$, $s\to-s$ and $t_0\to-t_0$.
Let $t_1:=(s^\alpha-t^\alpha+t_0^\alpha)^{1/\alpha}$,
then $s\geq t_1 \geq t_0$, because of the fact $s\geq t\geq t_0$.
From \eqref{xEEN3} it follows that for $s\geq t_1 \geq t_0\geq0$,
\begin{eqnarray*}
u((s^\alpha-t_1^\alpha+t_0^\alpha)^{1/\alpha})\leq \!\!\!\!\!\!\!\!\!&&N_2\frac{E_\alpha(\beta_2,t_0)}{E_\alpha(\beta_2,t_1)}\\
&&+\varepsilon N_1\int_{t_0}^{(s^\alpha-t_1^\alpha+t_0^\alpha)^{1/\alpha}}\tau^{\alpha-1}
\frac{E_\alpha(\beta_1,\tau)E_\alpha(\beta_1,t_1)}{E_\alpha(\beta_1,s)E_\alpha(\beta_1,t_0)}u(\tau)d\tau\\
&&+\varepsilon N_2\int_{(s^\alpha-t_1^\alpha+t_0^\alpha)^{1/\alpha}}^s\tau^{\alpha-1}
\frac{E_\alpha(\beta_2,s)E_\alpha(\beta_2,t_0)}{E_\alpha(\beta_2,\tau)E_\alpha(\beta_2,t_1)}u(\tau)d\tau.
\end{eqnarray*}
Put $v(t_1):=u((s^\alpha-t_1^\alpha+t_0^\alpha)^{1/\alpha})$ then $u(\tau)=v((s^\alpha-\tau^\alpha+t_0^\alpha)^{1/\alpha})$.
The inequality above yields that for $s\geq t_1 \geq t_0\geq0$,
\begin{eqnarray*}
v(t_1)\leq\!\!\!\!\!\!\!\!\!&& N_2\frac{E_\alpha(\beta_2,t_0)}{E_\alpha(\beta_2,t_1)}\\
&&+\varepsilon N_1\int_{t_0}^{(s^\alpha-t_1^\alpha+t_0^\alpha)^{1/\alpha}}\tau^{\alpha-1}
\frac{E_\alpha(\beta_1,\tau)E_\alpha(\beta_1,t_1)}{E_\alpha(\beta_1,s)E_\alpha(\beta_1,t_0)}v((s^\alpha-\tau^\alpha+t_0^\alpha)^{1/\alpha})d\tau\\
&&+\varepsilon N_2\int_{(s^\alpha-t_1^\alpha+t_0^\alpha)^{1/\alpha}}^s\tau^{\alpha-1}
\frac{E_\alpha(\beta_2,s)E_\alpha(\beta_2,t_0)}{E_\alpha(\beta_2,\tau)E_\alpha(\beta_2,t_1)}v((s^\alpha-\tau^\alpha+t_0^\alpha)^{1/\alpha})d\tau.
\end{eqnarray*}
Let $\iota:=(s^\alpha-\tau^\alpha+t_0^\alpha)^{1/\alpha}$, then
\begin{eqnarray*}
v(t_1)\leq\!\!\!\!\!\!\!\!\!&& N_2\frac{E_\alpha(\beta_2,t_0)}{E_\alpha(\beta_2,t_1)}
+\varepsilon N_2\int_{t_0}^{t_1}\iota^{\alpha-1}
\frac{E_\alpha(\beta_2,\iota)}{E_\alpha(\beta_2,t_1)}v(\iota)d\iota\\
&&+\varepsilon N_1\int_{t_1}^s\iota^{\alpha-1}
\frac{E_\alpha(\beta_1,t_1)}{E_\alpha(\beta_1,\iota)}v(\iota)d\iota,
~~~s\geq t_1 \geq t_0\geq0.
\end{eqnarray*}
The inequality also can be amplified as
\begin{eqnarray*}
v(t_1)\leq\!\!\!\!\!\!\!\!\!&& N_2\frac{E_\alpha(\beta_2,t_0)}{E_\alpha(\beta_2,t_1)}
+\frac{\varepsilon N_2}{E_\alpha(\beta_2,t_1)}{\cal I}^\alpha_{t_0}
E_\alpha(\beta_2,t_1)v(t_1)\\
&&-\varepsilon N_1 E_\alpha(\beta_1,t_1){\cal I}^\alpha_{+\infty}
\frac{v(t_1)}{E_\alpha(\beta_1,t_1)},
~~~t_1 \geq t_0\geq0.
\end{eqnarray*}
By the synchronous boundedness of both functions $u$ and $v$,
we employ Lemma \ref{efbds1} to gain
\begin{eqnarray*}
v(t_1)\leq K_2\frac{E_\alpha(\lambda_2,t_0)}{E_\alpha(\lambda_2,t_1)},~~~t_1\geq t_0.
\end{eqnarray*}
It follows from the definition of $v(t_1)$ that
\begin{eqnarray*}
u(t)\leq K_2\frac{E_\alpha(\lambda_2,t)}{E_\alpha(\lambda_2,s)},~~~s\geq t\geq t_0,
\end{eqnarray*}
where $K_2=\cfrac{N_2}{1-\theta}$, $\lambda_2=\beta_2-\cfrac{\varepsilon N_2}{1-\theta}$ given in Lemma \ref{efbds1}.
Hence, Corollary \ref{efbds2} is proved.
$\qquad\Box$

In the end of this section, we demonstrate the invariant manifolds theorem
for CFDE \eqref{TxAxft}.
But before we do that, let us introduce the following notion.

\begin{df}
\label{tan}
Let $\Omega$ is any subset of $\mathbb{R}^n$ including zero and $P$ is a projection matrix
such that $\mathbb{R}^n=P\mathbb{R}^n \oplus (I-P)\mathbb{R}^n$ and $P^2=P$.
We say $\Omega$ is tangent to $(I-P)\mathbb{R}^n$ (resp. $P\mathbb{R}^n$) at zero,
if $\|Px\|/\|(I-P)x\|\to 0$ (resp. $\|(I-P)x\|/\|Px\|\to 0$) as $x\to 0$ in $\Omega$.
\end{df}

From now on, let $k:=\mathcal{R}(I-P)$, where denote ${\cal R}(P)$ by the rank of matrix $P$
and assume that
\begin{description}
  \item[(E1)] $\zeta\in C_I(\mathbb{R}_+,\mathbb{R}_+)$ satisfies $\zeta(0)=0$;
  \item[(E2)] $\Lambda(\zeta)$ consists of functions $f\in C(\mathbb{R}^{n+1},\mathbb{R}^n)$ such that
    \begin{eqnarray*}
        f(t,0)\!\!\!\!\!\!\!\!&&=0,\\
        \|f(t,x)-f(t,y)\|\!\!\!\!\!\!\!\!&&\leq \zeta(\sigma)\|x-y\|,~~~\|x\|,\|y\|\leq\sigma;
    \end{eqnarray*}
  \item[(E3)] Projection matrix $P$ fulfils
  $X(t)P=PX(t)$ for all $t\in\mathbb{R}.$
\end{description}

\begin{thm}
Suppose that {\bf (E1)-(E3)} hold and
denote the unstable and stable manifolds of the hyperbolic equilibrium $x=0$ of equation \eqref{TxAxft}
as $U_k:=U_k(f)$ and $S_{n-k}:=S_{n-k}(f)$ respectively, for any $f\in \Lambda(\zeta)$.
Then $U_k$ and $S_{n-k}$ are tangent to $(I-P)\mathbb{R}^n$ and $P\mathbb{R}^n$ at $x=0$ respectively,
where $(I-P)\mathbb{R}^n$ and $P\mathbb{R}^n$ are the unstable and stable invariant subspaces
of the hyperbolic equilibrium $x=0$ of \eqref{TxAx}, respectively.
Moreover, there exist positive constants $M$, $\gamma_1$ and $\gamma_2$ such that
\begin{eqnarray}
\begin{split}
\|x(t)\|\leq &M\frac{E_\alpha(\gamma_1,t_0)}{E_\alpha(\gamma_1,t)}\|x(t_0)\|,~~~t\geq t_0\geq 0,~~~x(t_0)\in S_{n-k},\\
\|x(t)\|\leq &M\frac{E_\alpha(\gamma_2,t)}{E_\alpha(\gamma_2,t_0)}\|x(t_0)\|,~~~t\leq t_0\leq 0,~~~x(t_0)\in U_k.
\label{xMEE}
\end{split}
\end{eqnarray}
\label{xq}
\end{thm}

\begin{rmk}
The hyperbolic equilibrium $x=0$ of ODE $\dot{x}=A(t)x$ is also
the hyperbolic equilibrium of CFDE \eqref{TxAx}.
In fact, by Definition \ref{CFD}, one can verify
\begin{eqnarray*}
{\cal T}^\alpha x(0)=\lim_{t\to 0}{\cal T}^\alpha x(t)=\lim_{t\to 0}A(t)x(t)=\lim_{t\to 0}\dot{x}(t)=\dot{x}(0),
\end{eqnarray*}
which implies the assertion.
\end{rmk}

{\bf Proof of Theorem \ref{xq}.}
Assume $\lambda$, $N_i$, $\beta_i$ $(i=1,2)$ are all given in \eqref{XX-1KE} and \eqref{XPXK} respectively,
and the function $\zeta(\sigma)$ $(\sigma \geq 0)$ is given in {\bf (E1)}.
Take $\delta$ satisfy
\begin{eqnarray}
(\frac{N_1}{\beta_1}+\frac{N_2}{\beta_2})\zeta(\delta)<\frac{1}{2},
~~~N_1<(\beta_2+\lambda-4N_1N_2\zeta(\delta))(\frac{N_1}{\beta_1}+\frac{N_2}{\beta_2}).
\label{KK}
\end{eqnarray}
Choose $x_0$ satisfy $\|x_0\| \leq \delta/2N_1$ for $x_0\in \mathbb{R}^n$, and define
$\mathcal{L}(Px_0, \delta)$ is a set of functions $x\in C([t_0, +\infty),\mathbb{R}^n)$,
where $\|x\|_\infty:=\sup\limits_{t_0\leq t<+\infty}\|x(t)\|\leq \delta$ and $t_0\geq0$.
$\mathcal{L}(Px_0, \delta)$ is a closed bounded subset consisting of the Banach space
of all bounded continuous functions mapping $[t_0, +\infty)$ to $\mathbb{R}^n$ with the uniform topology.
For any $x\in \mathcal{L}(Px_0, \delta)$, define
\begin{eqnarray}
\begin{split}
({\cal J}x)(t):=&X(t)PX^{-1}(t_0)x_0+X(t)P{\cal I}^\alpha_{t_0}X^{-1}(t)f(t,x(t))\\
&+X(t)(I-P){\cal I}^\alpha_{+\infty}X^{-1}(t)f(t,x(t)),
~~~t\geq t_0.
\label{TxXPX}
\end{split}
\end{eqnarray}
It is easy to know that ${\cal J}x$ is well defined and continuous for $t\geq t_0$.
From \eqref{XPXK}, \eqref{KK} and {\bf (E2)}, we calculate that
\begin{eqnarray*}
\|({\cal J}x)(t)\|\leq \!\!\!\!\!\!\!\!&& N_1\frac{E_\alpha(\beta_1,t_0)}{E_\alpha(\beta_1,t)}\|x_0\|
+\frac{N_1}{E_\alpha(\beta_1,t)}{\cal I}^\alpha_{t_0}E_\alpha(\beta_1,t)\|f(t,x(t))\|
\\
&&-N_2E_\alpha(\beta_2,t){\cal I}^\alpha_{+\infty}\frac{\|f(t,x(t))\|}{E_\alpha(\beta_2,t)}
\\
\leq \!\!\!\!\!\!\!\!&& N_1\frac{E_\alpha(\beta_1,t_0)}{E_\alpha(\beta_1,t)}\|x_0\|
+\zeta(\delta)\Big{(}\frac{N_1}{\beta_1}+\frac{N_2}{\beta_2}\Big{)}\|x\|_\infty
\\
\leq \!\!\!\!\!\!\!\!&& N_1\|x_0\|+\zeta(\delta)\Big{(}\frac{N_1}{\beta_1}+\frac{N_2}{\beta_2}\Big{)}\delta
\\
< \!\!\!\!\!\!\!\!&& \frac{\delta}{2}+\frac{\delta}{2}=\delta,
\end{eqnarray*}
thus $\|{\cal J}x\|_\infty<\delta$ and ${\cal J}:\mathcal{L}(Px_0, \delta)\to \mathcal{L}(Px_0, \delta)$.

Analogously to the computation above, we obtain
\begin{eqnarray*}
\|({\cal J}x)(t)-({\cal J}y)(t)\|\leq \zeta(\delta)\Big{(}\frac{N_1}{\beta_1}+\frac{N_2}{\beta_2}\Big{)}\|x-y\|_\infty\leq \frac{1}{2}\|x-y\|_\infty,~~~t\geq t_0,
\end{eqnarray*}
which implies that ${\cal J}$ is a contraction mapping in $\mathcal{L}(Px_0, \delta)$.
In fact, there is a unique fixed point $x_*(t, Px_0)\in \mathcal{L}(Px_0, \delta)$ satisfying \eqref{xXPX1}.
Note that the function $x_*(t, Px_0)$ is continuous with respect to $Px_0$ and $x_*(t, 0)=0$.
Let $x_*(t):=x_*(t, Px_0)$ and $\hat{x}_*(t):=x_*(t, P\hat{x}_0)$, it follows from \eqref{xXPX1}, \eqref{XX-1KE} and {\bf (E3)} that
\begin{eqnarray*}
\|x_*(t)-\hat{x}_*(t)\|\leq \!\!\!\!\!\!\!\!&&
K\frac{E_\alpha(\lambda,t_0)}{E_\alpha(\lambda,t)}\|Px_0-P\hat{x}_0\|
\\
&&+\frac{N_1\zeta(\delta)}{E_\alpha(\beta_1,t)}{\cal I}^\alpha_{t_0} E_\alpha(\beta_1,t)\|x_*(t)-\hat{x}_*(t)\|
\\
&&-N_2\zeta(\delta)E_\alpha(\beta_2,t){\cal I}^\alpha_{+\infty}\frac{\|x_*(t)-\hat{x}_*(t)\|}{E_\alpha(\beta_2,t)},~~~t\geq t_0.
\end{eqnarray*}
By Lemma \ref{efbds1}, we acquire
\begin{eqnarray}
\|x_*(t, Px_0)-x_*(t, P\hat{x}_0)\|\leq 2K\cfrac{E_\alpha(\gamma_1,t_0)}{E_\alpha(\gamma_1,t)}\|Px_0-P\hat{x}_0\|,~~~t\geq t_0,
\label{xxK}
\end{eqnarray}
where $\gamma_1=\lambda-\cfrac{N_1\beta_1\beta_2}{N_1\beta_2+N_2\beta_1}$.
Combining the fact $x_*(t, 0)=0$ and \eqref{xxK},
one can verify that the first expression of the estimate \eqref{xMEE} is true.
Proceeding analogously to \eqref{xxK}, the second estimate in \eqref{xMEE} is also true,
when $\gamma_2=\lambda-\cfrac{N_2\beta_1\beta_2}{N_1\beta_2+N_2\beta_1}$.

Set $B_{\delta/2N_1}$ is the open ball in $\mathbb{R}^n$ centered at the origin with radius $\delta/2N_1$,
and take $S^*_{n-k}:=\{x:x=x_*(t_0,Px_0), x_0\in B_{\delta/2N_1}\cap\mathbb{R}^n\}$.
Let $h(Px_0):=x_*(t_0,Px_0)$ for $x_0\in B_{\delta/2N_1}\cap\mathbb{R}^n$,
we observe that $h$ is a continuous mapping from $B_{\delta/2N_1}\cap(P\mathbb{R}^n)$ to $S^*_{n-k}$,
then
\begin{eqnarray*}
h(Px_0)=Px_0+X(t_0)(I-P){\cal I}^\alpha_{+\infty}X^{-1}(t_0)f(t_0,x_*(t_0,Px_0)).
\end{eqnarray*}
Given $x_0,\hat{x}_0\in B_{\delta/2N_1}\cap\mathbb{R}^n$,
we employ \eqref{XPXK}, \eqref{KK}, \eqref{xxK} and {\bf (E2)} to obtain
\begin{eqnarray*}
\|h(Px_0)-h(P\hat{x}_0)\|\geq \!\!\!\!\!\!\!\!&&\|Px_0-P\hat{x}_0\|
+N_2\zeta(\delta)E_\alpha(\beta_2,t_0){\cal I}^\alpha_{+\infty}\frac{\|x_*(t_0)-\hat{x}_*(t_0)\|}{E_\alpha(\beta_2,t_0)}
\\
\geq \!\!\!\!\!\!\!\!&& \|Px_0-P\hat{x}_0\|\left[1+E_\alpha(\beta_2+\gamma_1,t_0){\cal I}^\alpha_{+\infty}\frac{2N_1N_2\zeta(\delta)}{E_\alpha(\beta_2+\gamma_1,t_0)}\right]
\\
\geq \!\!\!\!\!\!\!\!&&\left(1-\frac{2N_1N_2\zeta(\delta)}{\beta_2+\gamma_1}\right)\|Px_0-P\hat{x}_0\|
\geq\frac{1}{2}\|Px_0-P\hat{x}_0\|,
\end{eqnarray*}
yielding $h$ is a bijective.
And since $h^{-1}=P$ is continuous, $h$ is a homeomorphism.
Hence, $S^*_{n-k}$ is homeomorphic to the $(n-k)$-D open unit ball in $\mathbb{R}^{n-k}$.
If $S^*_{n-k}$ is not a positively invariant set,
then we expand $S^*_{n-k}$ into the positively invariant $S_{n-k}$,
by absorbing all the positive orbits of the solutions starting from $S^*_{n-k}$.
From the uniqueness of the solutions,
$S_{n-k}$ is also homeomorphic to the open unit ball in $\mathbb{R}^{n-k}$.
In other words, the case $\|Px\|<\delta/2N_1$ for all $x\in S_{n-k}$
implies $S_{n-k}\equiv S^*_{n-k}$.
It follows from \eqref{TxXPX}, \eqref{xxK}, {\bf (E2)} and the fact $x_*(t,0)=0$ that
\begin{eqnarray*}
\|(I-P)x_*(t_0,Px_0)\|\leq\!\!\!\!\!\!\!\!
&&-N_2 E_\alpha(\beta_2,t_0){\cal I}^\alpha_{+\infty}\frac{\|f(t_0,x_*(t_0,Px_0))\|}{E_\alpha(\beta_2,t_0)}
\\
\leq\!\!\!\!\!\!\!\!&& -N_2 E_\alpha(\beta_2,t_0){\cal I}^\alpha_{+\infty}\frac{\zeta(\|x_*(t_0,Px_0)\|)}{E_\alpha(\beta_2,t_0)}\|x_*(t_0,Px_0)\|
\\
\leq\!\!\!\!\!\!\!\!&& -N_2 E_\alpha(\beta_2,t_0){\cal I}^\alpha_{+\infty}\frac{\zeta(2N_1\|Px_0\|)}{E_\alpha(\beta_2,t_0)}2N_1\|Px_0\|
\\
\leq\!\!\!\!\!\!\!\!&& \frac{2N_1N_2}{\beta_2} \zeta(2N_1\|Px_0\|)\|Px_0\|.
\end{eqnarray*}
Since $\|Px_0\|\to 0$ as $\|x_0\|\to 0$, we get $\|(I-P)x_*(t_0,Px_0)\|/\|Px_0\|\to 0$ as $\|x_0\|\to 0$ in $S_{n-k}$.
Consequently $S_{n-k}$ is tangent to $P\mathbb{R}^n$ at zero. Similarly, one can construct the set $U_k$ via \eqref{xXPX2},
and complete the proof of Theorem \ref{xq}.
$\qquad\Box$

\section{Roughness of dichotomy}
\setcounter{equation}{0}

Our focus of this section is the roughness of the Mittag-Leffler dichotomy.
That is the preservation of dichotomy for hyperbolic linear systems
undergoing small linear perturbation.
Consider the perturbed equation of linear CFDE \eqref{TxAx} as follows
\begin{eqnarray}
{\cal T}^\alpha y=[A(t)+B(t)]y,~~~(t,y)\in\mathbb{R}^{n+1},
\label{TyABy}
\end{eqnarray}
where matrix functions $A\in C(\mathbb{R},\mathbb{R}^{n\times n})$
and $B\in C_b(\mathbb{R},\mathbb{R}^{n\times n})$.
The following is one of our main results of this paper.

\begin{thm}
Assume that $X(t)$ is a fundamental matrix of \eqref{TxAx} such that $X(0)=I$,
and equation \eqref{TxAx} possesses a Mittag-Leffler dichotomy,
i.e., estimates \eqref{XPXK} hold in $\mathbb{R}_+$.
If $\varepsilon:=\sup\limits_{t\geq 0}\|B(t)\|$ is sufficiently small,
then perturbed equation \eqref{TyABy} also possesses a Mittag-Leffler dichotomy in $\mathbb{R}_+$.
\label{ccx}
\end{thm}

{\bf Proof.}
We divide the proof of Theorem \ref{ccx} into the following three steps.

Step 1: Finding bounded solutions of equation \eqref{TyABy}.
Let matrix function $Y\in C_b(\mathbb{R}_+,\mathbb{R}^{n\times n})$ equipped with norm
\begin{eqnarray*}
\|Y\|_\infty:=\sup\limits_{t\geq 0}\|Y(t)\|.
\end{eqnarray*}
Define mapping $L:C_b(\mathbb{R}_+,\mathbb{R}^{n\times n})\to C_b(\mathbb{R}_+,\mathbb{R}^{n\times n})$
as
\begin{eqnarray*}
LY(t)=X(t)P+X(t)P{\cal I}^\alpha_{0}X^{-1}(t)B(t)Y(t)+X(t)(I-P){\cal I}^\alpha_{+\infty}X^{-1}(t)B(t)Y(t).
\end{eqnarray*}
It follows from \eqref{XPXK} that
\begin{eqnarray*}
\|LY(t)\|\leq\frac{N_1}{E_\alpha(\beta_1,t)}
+\frac{\varepsilon N_1\|Y\|_\infty}{E_\alpha(\beta_1,t)}{\cal I}^\alpha_{0}E_\alpha(\beta_1,t)
-E_\alpha(\beta_2,t){\cal I}^\alpha_{+\infty}\frac{\varepsilon N_2\|Y\|_\infty}{E_\alpha(\beta_2,t)}.
\end{eqnarray*}
Observing that $LY(t)$ is bounded and continuous for $t\geq0$, we obtain
\begin{eqnarray*}
\|LY\|_\infty\leq N_1+\varepsilon\Big{(}\frac{N_1}{\beta_1}+\frac{N_2}{\beta_2}\Big{)}\|Y\|_\infty.
\end{eqnarray*}
Given another $\hat{Y}\in C_b(\mathbb{R}_+,\mathbb{R}^{n\times n})$,
analogously we get
\begin{eqnarray*}
\|LY-L\hat{Y}\|_\infty\leq \varepsilon\Big{(}\frac{N_1}{\beta_1}+\frac{N_2}{\beta_2}\Big{)}\|Y-\hat{Y}\|_\infty.
\end{eqnarray*}
This yields that the mapping $L$ has a unique $Y_1\in C_b(\mathbb{R}_+,\mathbb{R}^{n\times n})$
such that
\begin{eqnarray}
\begin{split}
Y_1(t)=&X(t)P+X(t)P{\cal I}^\alpha_{0}X^{-1}(t)B(t)Y_1(t)\\
&+X(t)(I-P){\cal I}^\alpha_{+\infty}X^{-1}(t)B(t)Y_1(t),
\end{split}
\label{y1t}
\end{eqnarray}
if
\begin{eqnarray*}
\theta:=\varepsilon\Big{(}\frac{N_1}{\beta_1}+\frac{N_2}{\beta_2}\Big{)}<1.
\end{eqnarray*}
Obviously, $Y_1(t)$ is also a matrix solution of \eqref{TyABy} and differentiable.
Post-projecting $P$ on both hands sides of \eqref{y1t},
we also know that $Y_1(t)P$ is the unique fixed point of $L$,
and $Y_1(t)P=Y_1(t)$.

Step 2: Constructing projection matrix.
Let $Q:=Y_1(0)$, then $QP=Q$.
Combining \eqref{y1t} with the property $P(I-P)=0$, and replacing $t$ with $s$, we attain
\begin{eqnarray}
X(t)PX^{-1}(s)Y_1(s)=X(t)P+X(t)P{\cal I}^\alpha_{0}X^{-1}(s)B(s)Y_1(s).
\label{xpxy}
\end{eqnarray}
It follows from \eqref{y1t} and \eqref{xpxy} that
\begin{eqnarray}
\begin{split}
Y_1(t)=&X(t)PX^{-1}(s)Y_1(s)+X(t)P{\cal I}^\alpha_{s}X^{-1}(t)B(t)Y_1(t)\\
&+X(t)(I-P){\cal I}^\alpha_{+\infty}X^{-1}(t)B(t)Y_1(t),
~~~t\geq s\geq0.
\end{split}
\label{YtXPX}
\end{eqnarray}
Noting \eqref{xpxy} with $t=s=0$, we gain $PQ=P$.
Post-projecting $Q$ on both hands sides of \eqref{y1t} again,
we acquire
\begin{eqnarray*}
Y_1(t)Q=\!\!\!\!\!\!\!\!&&X(t)P+X(t)P{\cal I}^\alpha_{0}X^{-1}(t)B(t)Y_1(t)Q
\\
&&+X(t)(I-P){\cal I}^\alpha_{+\infty}X^{-1}(t)B(t)Y_1(t)Q,
\end{eqnarray*}
implying $Y_1(t)Q$ is also a fixed point of $L$. In conclusion,
\begin{eqnarray*}
Y_1(t)Q=Y_1(t)=Y_1(t)P.
\end{eqnarray*}
Obviously, $Q$ is a projection when $t=0$.

Provided $Y(t)$ is a fundamental matrix of \eqref{TyABy} fulfilling $Y(0)=I$, we derive
\begin{eqnarray}
Y_1(t)=Y(t)Q.
\label{Y1YQ}
\end{eqnarray}
Set
\begin{eqnarray}
Y_2(t):=Y(t)(I-Q),
\label{Y2YIQ}
\end{eqnarray}
then $Y(t)=Y_1(t)+Y_2(t)$.
Relying on the variation of constants formula \eqref{csbygs}, we calculate that
\begin{eqnarray}
Y_2(t)\!\!\!\!\!\!\!\!&&=X(t)X^{-1}(0)Y_2(0)+X(t){\cal I}^\alpha_{0}X^{-1}(t)B(t)Y_2(t)
\nonumber\\
\!\!\!\!\!\!\!\!&&=X(t)Y(0)(I-Q)+X(t){\cal I}^\alpha_{0}X^{-1}(t)B(t)Y_2(t)
\nonumber\\
\!\!\!\!\!\!\!\!&&=X(t)(I-Q)+X(t){\cal I}^\alpha_{0}X^{-1}(t)B(t)Y_2(t).
\label{Y2XIQ}
\end{eqnarray}
Combining \eqref{Y2XIQ} with the fact $(I-P)(I-Q)=I-Q$, and replacing $t$ with $s$, we acquire
\begin{eqnarray}
X(t)(I-P)X^{-1}(s)Y_2(s)=X(t)(I-Q)+X(t)(I-P){\cal I}^\alpha_{0}X^{-1}(s)B(s)Y_2(s).
\label{XIPXY2}
\end{eqnarray}
Subsequently, by \eqref{Y2XIQ} and \eqref{XIPXY2}, we receive
\begin{eqnarray}
Y_2(t)=\!\!\!\!\!\!\!\!&&X(t)(I-P)X^{-1}(s)Y_2(s)+X(t){\cal I}^\alpha_{0}X^{-1}(t)B(t)Y_2(t)
\nonumber\\
&&-X(t)(I-P){\cal I}^\alpha_{0}X^{-1}(s)B(s)Y_2(s)
\nonumber\\
=\!\!\!\!\!\!\!\!&&X(t)(I-P)X^{-1}(s)Y_2(s)+X(t)P{\cal I}^\alpha_{0}X^{-1}(t)B(t)Y_2(t)
\nonumber\\
&&+X(t)(I-P){\cal I}^\alpha_{0}X^{-1}(t)B(t)Y_2(t)-X(t)(I-P){\cal I}^\alpha_{0}X^{-1}(s)B(s)Y_2(s)
\nonumber\\
=\!\!\!\!\!\!\!\!&&X(t)(I-P)X^{-1}(s)Y_2(s)+X(t)P{\cal I}^\alpha_{0}X^{-1}(t)B(t)Y_2(t)
\nonumber\\
&&+X(t)(I-P){\cal I}^\alpha_{s}X^{-1}(t)B(t)Y_2(t),~~~s\geq t\geq 0.
\label{Y2XIPX}
\end{eqnarray}
From \eqref{YtXPX} and \eqref{Y2XIPX} it follows that for any vector $\xi$,
\begin{eqnarray*}
\|Y_1(t)\xi\|\leq\!\!\!\!\!\!\!\!&& N_1\frac{E_\alpha(\beta_1,s)}{E_\alpha(\beta_1,t)}\|Y_1(s)\xi\|
+\frac{\varepsilon N_1}{E_\alpha(\beta_1,t)}|{\cal I}^\alpha_{s}E_\alpha(\beta_1,t)Y_1(t)\xi|
\\
&&-\varepsilon N_2E_\alpha(\beta_2,t)\Big{|}{\cal I}^\alpha_{+\infty}\frac{Y_1(t)\xi}{E_\alpha(\beta_2,t)}\Big{|},
~~~t\geq s\geq0,
\end{eqnarray*}
and
\begin{eqnarray*}
\|Y_2(t)\xi\|\leq\!\!\!\!\!\!\!\!&& N_2\frac{E_\alpha(\beta_2,t)}{E_\alpha(\beta_2,s)}\|Y_2(s)\xi\|
+\frac{\varepsilon N_1}{E_\alpha(\beta_1,t)}|{\cal I}^\alpha_{0}E_\alpha(\beta_1,t)Y_2(t)\xi|
\\
&&-\varepsilon N_2E_\alpha(\beta_2,t)\Big{|}{\cal I}^\alpha_{s}\frac{Y_2(t)\xi}{E_\alpha(\beta_2,t)}\Big{|},
~~~s\geq t\geq0.
\end{eqnarray*}
Thus, by Lemma \ref{efbds1} and Corollary \ref{efbds2}, we can know
\begin{eqnarray}
\begin{split}
\|Y_1(t)\xi\|\leq K_1\frac{E_\alpha(\lambda_1,s)}{E_\alpha(\lambda_1,t)}\|Y_1(s)\xi\|,
~~~t\geq s\geq0,
\\
\|Y_2(t)\xi\|\leq K_2\frac{E_\alpha(\lambda_2,t)}{E_\alpha(\lambda_2,s)}\|Y_2(s)\xi\|,
~~~s\geq t\geq0,
\label{Y1Y2}
\end{split}
\end{eqnarray}
where
$K_i=\cfrac{N_i}{1-\theta}$, $\lambda_i=\beta_i-\cfrac{\varepsilon N_i}{1-\theta}$ and $\theta=\varepsilon\Big{(}\cfrac{N_1}{\beta_1}+\cfrac{N_2}{\beta_2}\Big{)}$
for $i=1,2$.

Step 3: Estimation of fundamental solutions.
To prove from \eqref{Y1Y2} that the perturbed equation \eqref{TyABy} also possesses a Mittag-Leffler dichotomy,
we only need to exhibit that $Y(t)QY^{-1}(t)$ is bounded.
From the facts $(I-P)P=0$, $(I-P)(I-P)=I-P$ and \eqref{y1t} it follows that
\begin{eqnarray*}
X(t)(I-P)X^{-1}(t)Y_1(t)=X(t)(I-P){\cal I}^\alpha_{+\infty}X^{-1}(t)B(t)Y_1(t).
\end{eqnarray*}
By \eqref{Y1Y2}, for any vector $\xi$, we calculate that
\begin{eqnarray}
\|X(t)(I-P)X^{-1}(t)Y_1(t)\xi\|\!\!\!\!\!\!\!\!&&\leq -\varepsilon N_2 E_\alpha(\beta_2,t){\cal I}^\alpha_{+\infty}\frac{\|Y_1(t)\xi\|}{E_\alpha(\beta_2,t)}
\nonumber\\
&&\leq -E_\alpha(\lambda_1+\beta_2,t)\|Y_1(t)\xi\|{\cal I}^\alpha_{+\infty}\frac{\varepsilon K_1 N_2}{E_\alpha(\lambda_1+\beta_2,t)}
\nonumber\\
&&\leq \frac{\varepsilon K_1 N_2}{\lambda_1+\beta_2}\|Y_1(t)\xi\|.
\label{XIPXY1}
\end{eqnarray}
Analogously, pre-multiplying $X(t)PX^{-1}(t)$ on both hands sides of \eqref{Y2XIQ},
by the property $P(I-Q)=0$, we obtain
\begin{eqnarray*}
X(t)PX^{-1}(t)Y_2(t)=X(t)P{\cal I}^\alpha_{0}X^{-1}(t)B(t)Y_2(t).
\end{eqnarray*}
It follows from \eqref{Y1Y2} that for any vector $\xi$,
\begin{eqnarray}
\|X(t)PX^{-1}(t)Y_2(t)\xi\|\!\!\!\!\!\!\!\!&&\leq \frac{\varepsilon N_1}{E_\alpha(\beta_1,t)}\Big{|}{\cal I}^\alpha_{0}E_\alpha(\beta_1,t)Y_2(t)\xi\Big{|}
\nonumber\\
&&\leq\frac{\varepsilon K_2 N_1}{E_\alpha(\lambda_2+\beta_1,t)}\|Y_2(t)\xi\|{\cal I}^\alpha_{0}E_\alpha(\lambda_2+\beta_1,t)
\nonumber\\
&&\leq\frac{\varepsilon K_2 N_1}{\lambda_2+\beta_1}\|Y_2(t)\xi\|.
\label{XPXY2}
\end{eqnarray}
Substituting \eqref{Y1YQ} and \eqref{Y2YIQ} into \eqref{XIPXY1} and \eqref{XPXY2} respectively,
and replacing $\xi$ by $Y^{-1}(t)\xi$, we acquire
\begin{eqnarray}
\|X(t)(I-P)X^{-1}(t)Y_1(t)\xi\|\!\!\!\!\!\!\!\!&&\leq\|X(t)(I-P)X^{-1}(t)Y(t)QY^{-1}(t)\xi\|
\nonumber\\
&&\leq\frac{\varepsilon K_1 N_2}{\lambda_1+\beta_2}\|Y(t)QY^{-1}(t)\xi\|,
\label{XI-PX}
\end{eqnarray}
and
\begin{eqnarray}
\|X(t)PX^{-1}(t)Y_2(t)\xi\|\!\!\!\!\!\!\!\!&&\leq\|X(t)PX^{-1}(t)Y(t)(I-Q)Y^{-1}(t)\xi\|
\nonumber\\
&&\leq\frac{\varepsilon K_2 N_1}{\lambda_2+\beta_1}\|Y(t)(I-Q)Y^{-1}(t)\xi\|.
\label{XPX-1}
\end{eqnarray}

On the other hand, it is evident to derive that
\begin{eqnarray}
Y(t)QY^{-1}(t)-X(t)PX^{-1}(t)=\!\!\!\!\!\!\!\!&&X(t)[P+(I-P)]X^{-1}(t)Y(t)QY^{-1}(t)
\nonumber\\
&&-X(t)PX^{-1}(t)Y(t)[Q+(I-Q)]Y^{-1}(t)
\nonumber\\
=\!\!\!\!\!\!\!\!&&X(t)PX^{-1}(t)Y(t)QY^{-1}(t)
\nonumber\\
&&+X(t)(I-P)X^{-1}(t)Y(t)QY^{-1}(t)
\nonumber\\
&&-X(t)PX^{-1}(t)Y(t)QY^{-1}(t)
\nonumber\\
&&-X(t)PX^{-1}(t)Y(t)(I-Q)Y^{-1}(t)
\nonumber\\
=\!\!\!\!\!\!\!\!&&X(t)(I-P)X^{-1}(t)Y(t)QY^{-1}(t)
\nonumber\\
&&-X(t)PX^{-1}(t)Y(t)(I-Q)Y^{-1}(t).
\label{YQYXPX}
\end{eqnarray}
Combining \eqref{XI-PX} and \eqref{XPX-1} with \eqref{YQYXPX},
we can obtain
\begin{eqnarray*}
\|Y(t)QY^{-1}(t)-X(t)PX^{-1}(t)\|\leq \frac{\varepsilon K_1 N_2}{\lambda_1+\beta_2}\mu_1+\frac{\varepsilon K_2 N_1}{\lambda_2+\beta_1}\mu_2,
\end{eqnarray*}
where $\mu_1(t):=\|Y(t)QY^{-1}(t)\|$, $\mu_2(t):=\|Y(t)(I-Q)Y^{-1}(t)\|$. For convenience,
take $N:=\max\{N_1,N_2\}$ and $\beta:=\min\{\beta_1,\beta_2\}$ such that $\theta\leq \hat{\theta}:=2\varepsilon N/\beta$,
and it yields
\begin{eqnarray}
\mu_1=\|Y(t)QY^{-1}(t)\|\!\!\!\!\!\!\!\!&&\leq \Big{(}\frac{\varepsilon K_1 N_2}{\lambda_1+\beta_2}\mu_1
+\frac{\varepsilon K_2 N_1}{\lambda_2+\beta_1}\mu_2\Big{)}+\|X(t)PX^{-1}(t)\|
\nonumber\\
&&\leq \Big{(}\frac{\varepsilon K_1 N_2}{\lambda_1+\beta_2}\mu_1
+\frac{\varepsilon K_2 N_1}{\lambda_2+\beta_1}\mu_2\Big{)}+N_1
\nonumber\\
&&\leq \eta(\mu_1+\mu_2)+N,
\label{lamda1}
\end{eqnarray}
where $\eta=\cfrac{\varepsilon N^2}{2\beta-5\varepsilon N}$.
It is obvious that
\begin{eqnarray*}
Y(t)QY^{-1}(t)-X(t)PX^{-1}(t)=X(t)(I-P)X^{-1}(t)-Y(t)(I-Q)Y^{-1}(t),
\end{eqnarray*}
then
\begin{eqnarray}
\mu_2\!\!\!\!\!\!\!\!&&\leq \Big{(}\frac{\varepsilon K_1 N_2}{\lambda_1+\beta_2}\mu_1
+\frac{\varepsilon K_2 N_1}{\lambda_2+\beta_1}\mu_2\Big{)}+N_2
\leq \eta(\mu_1+\mu_2)+N.
\label{lamda2}
\end{eqnarray}
By adding the inequality \eqref{lamda1} and \eqref{lamda2}, we attain
\begin{eqnarray*}
\mu_1+\mu_2\leq \frac{2N}{1-2\eta}.
\end{eqnarray*}
If $\eta<1/2$, then
\begin{eqnarray*}
\mu_1,\mu_2\leq \eta(\mu_1+\mu_2)+N
\leq \frac{N}{1-2\eta}.
\end{eqnarray*}
Substituting \eqref{Y1YQ} and \eqref{Y2YIQ} into \eqref{Y1Y2},
and replacing $\xi$ by $Y^{-1}(s)\xi$, we gain
\begin{eqnarray*}
\|Y(t)QY^{-1}(s)\xi\|\!\!\!\!\!\!\!\!&&\leq K_1\frac{E_\alpha(\lambda_1,s)}{E_\alpha(\lambda_1,t)}\|Y(s)QY^{-1}(s)\xi\|
\\
&&\leq \frac{K_1N}{1-2\eta}\frac{E_\alpha(\lambda_1,s)}{E_\alpha(\lambda_1,t)}\|\xi\|,
~~~t\geq s\geq0,
\\
\|Y(t)(I-Q)Y^{-1}(s)\xi\|\!\!\!\!\!\!\!\!&&\leq  K_2\frac{E_\alpha(\lambda_2,t)}{E_\alpha(\lambda_2,s)}\|Y(s)(I-Q)Y^{-1}(s)\xi\|
\\
&&\leq \frac{K_2N}{1-2\eta}\frac{E_\alpha(\lambda_2,t)}{E_\alpha(\lambda_2,s)}\|\xi\|,
~~~s\geq t\geq0.
\end{eqnarray*}
For the arbitrariness of vector $\xi$, we obtain the Mittag-Leffler dichotomy as follows
\begin{eqnarray*}
\|Y(t)QY^{-1}(s)\|\leq \frac{K_1N}{1-2\eta}\frac{E_\alpha(\lambda_1,s)}{E_\alpha(\lambda_1,t)},
~~~t\geq s\geq0,
\\
\|Y(t)(I-Q)Y^{-1}(s)\|\leq \frac{K_2N}{1-2\eta}\frac{E_\alpha(\lambda_2,t)}{E_\alpha(\lambda_2,s)},
~~~s\geq t\geq0.
\end{eqnarray*}
Therefore, Theorem \ref{ccx} is proved completely.
$\qquad\Box$

Finally, we present the concrete constants of estimates in Theorem \ref{ccx}.
Like the condition $N=\max\{N_1,N_2\}$ and $\beta=\min\{\beta_1,\beta_2\}$
such that $\theta\leq \hat{\theta}=2\varepsilon N/\beta$ holds,
let $N\geq 1$ and $\hat{\theta}<2/5N$ such that
$\eta=\cfrac{\varepsilon N^2}{2\beta-5\varepsilon N}<\cfrac{N}{10N-5}<\cfrac{1}{2}$.
Thus, by elementary calculation, we can obtain the following brief statement.

\begin{co}
Suppose that equation \eqref{TxAx} possesses the Mittag-Leffler dichotomy \eqref{XPXK} in $\mathbb{R}_+$.
If
\begin{eqnarray*}
\varepsilon:=\sup\limits_{t\in \mathbb{R}_+}\|B(t)\|<\frac{\beta}{5N^2},
\end{eqnarray*}
then perturbed equation \eqref{TyABy} also possesses the following Mittag-Leffler dichotomy:
\begin{eqnarray*}
\|Y(t)QY^{-1}(s)\|\leq \frac{25N^2}{9} \frac{E_\alpha(\beta-3\varepsilon N,s)}{E_\alpha(\beta-3\varepsilon N,t)},~~~t\geq s\geq 0,\\
\|Y(t)(I-Q)Y^{-1}(s)\|\leq \frac{25N^2}{9} \frac{E_\alpha(\beta-3\varepsilon N,t)}{E_\alpha(\beta-3\varepsilon N,s)},~~~s\geq t \geq 0,
\end{eqnarray*}
where $Y(t)$ is a fundamental matrix of \eqref{TyABy} such that $Y(0)=I$, and both projection matrices $Q$ and $P$ have the same rank.
Moreover,
\begin{eqnarray*}
\|Y(t)QY^{-1}(t)-X(t)PX^{-1}(t)\|\leq \frac{2\varepsilon N^3}{2\beta-5\varepsilon N-2\varepsilon N^2},
~~~t\geq 0.
\end{eqnarray*}
\end{co}

\section{Nonuniform dichotomy}
\setcounter{equation}{0}
This section is a continuation of studies for the Mittag-Leffler dichotomy.
More precisely, we concern {\it nonuniform Mittag-Leffler dichotomy}.
Let $\mathcal{B}(Z)$ consist of all bounded linear operators in Banach space $Z$.
Consider nonautonomous linear CFDE on $Z$
\begin{eqnarray}
{\cal T}^\alpha x=A(t)x,~~~(t,x)\in J\times Z,
\label{TxAxj}
\end{eqnarray}
where linear operator $A\in C(J,\mathcal{B}(Z))$ for some interval $J\subset \mathbb{R}$
and $\mathcal{B}(Z)$ is also a Banach space with the norm
$\|A\|:=\sup\limits_{x\in Z,\|x\|=1}\|Ax\|$ for all $A\in \mathcal{B}(Z)$.
Let $T(t,s)$ be a family of evolution operators satisfying $x(t)=T(t,s)x(s)$ for $t\geq s$ and $t,s\in J$,
where $x(t)$ is any solution of \eqref{TxAxj}. $T(t,s)$ further satisfies:
\begin{description}
  \item[(F1)] $T(t,t)=\mathrm{Id}$ (abbreviation of identity) for $t\in J$;
  \item[(F2)] $T(t,s)T(s,\tau)=T(t,\tau)$ for $t,s,\tau\in J$;
  \item[(F3)] the evolution operator $T(t,s)$ is invertible and $T^{-1}(t,s)=T(s,t)$ for $t,s\in J$.
\end{description}

First, we introduce the notions of nonuniform asymptotical stability
and nonuniform Mittag-Leffler dichotomy.

\begin{df}
Equation \eqref{TxAxj} is said to be nonuniform asymptotically stable in $J$
if there exist constants $\hat{N},\hat{\beta}>0$ and $\epsilon \geq 0$ such that
\begin{eqnarray}
\|T(t,s)\|\leq \hat{N}\frac{E_\alpha(\hat{\beta},s)}{E_\alpha(\hat{\beta},t)}E_\alpha(\epsilon,|s|),~~~t\geq s,~~~t,s\in J.
\label{TKE}
\end{eqnarray}
In particular, \eqref{TxAxj} is uniformly asymptotically stable like \eqref{XX-1KE} if \eqref{TKE} hold with $\epsilon=0$.
\end{df}
\begin{df}
Equation \eqref{TxAxj} is said to admit a nonuniform Mittag-Leffler dichotomy in $J$
if there exist projections $P:J\to \mathcal{B}(Z)$ such that
\begin{eqnarray}
T(t,s)P(s)=P(t)T(t,s),~~~t\geq s,~~~t,s\in J,
\label{TPPT}
\end{eqnarray}
and constants $\hat{N}_i,\hat{\beta}_i>0$ $(i=1,2)$ and $\epsilon \geq 0$ such that for $t,s\in J$,
\begin{eqnarray}
\begin{split}
\|T(t,s)P(s)\|\leq \hat{N}_1\frac{E_\alpha(\hat{\beta}_1,s)}{E_\alpha(\hat{\beta}_1,t)}E_\alpha(\epsilon,|s|),~~~t\geq s,
\\
\|T(t,s)(\mathrm{Id}-P(s))\|\leq \hat{N}_2\frac{E_\alpha(\hat{\beta}_2,t)}{E_\alpha(\hat{\beta}_2,s)}E_\alpha(\epsilon,|s|),~~~s\geq t.
\label{TPK}
\end{split}
\end{eqnarray}
In particular, \eqref{TxAxj} admits a uniform Mittag-Leffler dichotomy like \eqref{XPXK} if \eqref{TPK} hold with $\epsilon=0$.
\label{fyz}
\end{df}

All results in this section are presented in $\mathbb{R}_+$,
and denote
\begin{eqnarray*}
{\cal I}^\alpha_{s}f(t,\cdot):=\int_s^t\tau^{\alpha-1}f(t,\tau)d\tau.
\end{eqnarray*}
Consider the linear perturbation of \eqref{TxAxj} as follows
\begin{eqnarray}
{\cal T}^\alpha x=[A(t)+B(t)]x,~~~(t,x)\in \mathbb{R}_+\times Z,
\label{TxABx}
\end{eqnarray}
where linear operators $A\in C(\mathbb{R}_+,\mathcal{B}(Z))$
and $B\in C_b(\mathbb{R}_+,\mathcal{B}(Z))$.
The following theorem gives out the roughness of nonuniform asymptotical stability.

\begin{thm}
Assume that equation \eqref{TxAxj} admits nonuniform asymptotical stability in $\mathbb{R}_+$,
and there exists constant $\delta$ such that $\|B(t)\|\leq \delta/E_\alpha(\epsilon,t)$ for $t\in \mathbb{R}_+$.
If $\theta:=\delta\hat{N}/\hat{\beta}<1$, then equation \eqref{TxABx} also admits nonuniform asymptotical stability in $\mathbb{R}_+$,
that is,
\begin{eqnarray}
\|U(t,s)\|\leq \frac{\hat{N}}{1-\theta} \frac{E_\alpha(\gamma,s)}{E_\alpha(\gamma,t)}E_\alpha(\epsilon,s),~~~t\geq s,~~~t,s\in \mathbb{R}_+,
\label{UKEE}
\end{eqnarray}
where $\gamma=\hat{\beta}-\cfrac{\delta \hat{N}}{1-\theta}$ and $U(t,s)$ denotes the evolution operator associated to \eqref{TxABx}.
\label{fyzjjys}
\end{thm}

{\bf Proof.}
Consider the space
\begin{center}
$W:=\{U(t,s)_{t\geq s}\in \mathcal{B}(Z):U$ is continuous and $\|U\|_\alpha<\infty$, $(t,s)\in \mathbb{R}^2_+$$\}$,
\end{center}
equipped with $\alpha$-weighted norm
\begin{eqnarray}
\|U\|_\alpha:=\sup\left\{\frac{\|U(t,s)\|}{E_\alpha(\epsilon,s)}:t\geq s,(t,s)\in \mathbb{R}^2_+\right\}.
\label{uslts}
\end{eqnarray}
It is easy to verify that $W$ is a Banach space.
In $W$ define an operator ${\cal J}$ by
\begin{eqnarray*}
({\cal J}U)(t,s)=T(t,s)+{\cal I}^\alpha_{s}T(t,\cdot)B(\cdot)U(\cdot,s).
\end{eqnarray*}
It follows from \eqref{TKE} that
\begin{eqnarray*}
\|({\cal J}U)(t,s)\|\!\!\!\!\!\!\!\!&&\leq\|T(t,s)\|+{\cal I}^\alpha_{s}\|T(t,\cdot)\| \|B(\cdot)\| \|U(\cdot,s)\|
\\
&&\leq \hat{N}\frac{E_\alpha(\hat{\beta},s)}{E_\alpha(\hat{\beta},t)}E_\alpha(\epsilon,s)
+\frac{\delta \hat{N} \|U\|_\alpha E_\alpha(\epsilon,s)}{E_\alpha(\hat{\beta},t)}{\cal I}^\alpha_{s}E_\alpha(\hat{\beta},t)
\\
&&\leq \hat{N}E_\alpha(\epsilon,s)
+\frac{\delta\hat{N}}{\hat{\beta}} \|U\|_\alpha E_\alpha(\epsilon,s).
\end{eqnarray*}
And by \eqref{uslts} we obtain
\begin{eqnarray*}
\|{\cal J}U\|_\alpha\leq \hat{N}+\frac{\delta \hat{N}}{\hat{\beta}}\|U\|_\alpha<\infty,
\end{eqnarray*}
which yields that the operator ${\cal J}:W\to W$ is well defined.
Analogously to the computation above, we have
\begin{eqnarray*}
\|{\cal J}U_1-{\cal J}U_2\|_\alpha\leq \frac{\delta \hat{N}}{\hat{\beta}}\|U_1-U_2\|_\alpha,~~~U_1,U_2\in W,
\end{eqnarray*}
which implies that ${\cal J}$ is a contraction since $\delta<\hat{\beta}/\hat{N}$.
So there exists a unique $U\in W$ satisfying ${\cal J}U=U$,
and one can verify that it is a solution of \eqref{TxABx}.
We apply Lemma \ref{efbds1} with condition $\theta:=\delta \hat{N}/\hat{\beta}<1$
to the estimation of $\|U(t,s)\|$.
And inequality \eqref{UKEE} is true.
$\qquad\Box$

Subsequently,
our purpose is to establish roughness of nonuniform Mittag-Leffler dichotomy in $\mathbb{R}_+$.
A preliminary theorem and the main theorem of roughness are both stated as follows.

\begin{thm}
Assume that equation \eqref{TxAxj} admits a nonuniform Mittag-Leffler dichotomy \eqref{TPK} in $\mathbb{R}_+$,
and there exists constant $\delta$ such that $\|B(t)\|\leq \delta/E_\alpha(\epsilon,t)$ for $t\in \mathbb{R}_+$.
If
\begin{eqnarray}
\theta:=\delta\left(\frac{\hat{N}_1}{\hat{\beta}_1}+\frac{\hat{N}_2}{\hat{\beta}_2}\right)<1,~~~\epsilon< \min\{\hat{\beta}_1,\hat{\beta}_2\},
\label{the}
\end{eqnarray}
then there exist projections $\hat{P}:\mathbb{R}_+\to \mathcal{B}(Z)$ such that
\begin{eqnarray}
\hat{T}(t,s)\hat{P}(s)=\hat{P}(t)\hat{T}(t,s),~~~t\geq s,~~~t,s\in \mathbb{R}_+,
\label{Tts}
\end{eqnarray}
and constants $K_i,\lambda_i>0$ $(i=1,2)$ and $\epsilon\geq 0$ such that
\begin{eqnarray}
\begin{split}
\|\hat{T}(t,s)|\mathrm{Im}\hat{P}(s)\|\leq K_1\frac{E_\alpha(\lambda_1,s)}{E_\alpha(\lambda_1,t)}E_\alpha(\epsilon,s),~~~t\geq s\geq 0,
\\
\|\hat{T}(t,s)|\mathrm{Im}(\mathrm{Id}-\hat{P}(s))\|\leq K_2\frac{E_\alpha(\lambda_2,t)}{E_\alpha(\lambda_2,s)}E_\alpha(\epsilon,s),~~~s\geq t\geq 0,
\label{TPKE}
\end{split}
\end{eqnarray}
where $K_i=\cfrac{\hat{N}_i}{1-\theta}$, $\lambda_i=\hat{\beta}_i-\cfrac{\delta \hat{N}_i}{1-\theta}$ $(i=1,2)$,
and $\hat{T}(t,s)$ is the evolution operator associated to equation \eqref{TxABx}.
\label{thmfyz}
\end{thm}

\begin{thm}
Assume that equation \eqref{TxAxj} admits a nonuniform Mittag-Leffler dichotomy \eqref{TPK} in $\mathbb{R}_+$
under condition \eqref{the}.
If $\delta$ is sufficiently small such that $\|B(t)\|\leq \delta/E_\alpha(2\epsilon,t)$ for $t\in \mathbb{R}_+$,
then equation \eqref{TxABx} also  admits a nonuniform Mittag-Leffler dichotomy in $\mathbb{R}_+$.
\label{thmfyz1}
\end{thm}

{\bf Proof of Theorem \ref{thmfyz}.}
We divide the proof into the following several steps.

Step 1: Construction of bounded solutions for \eqref{TxABx}.
Recall space $W$ in Theorem \ref{fyzjjys},
then the following lemma gives out the existence of bounded solution.

\begin{lm}
For each $t,s\in \mathbb{R}_+$, equation \eqref{TxABx} has a unique solution $U\in W$ such that
\begin{eqnarray}
\begin{split}
U(t,s)=&T(t,s)P(s)+{\cal I}^\alpha_{s}T(t,\cdot)P(\cdot)B(\cdot)U(\cdot,s)
\\
&+{\cal I}^\alpha_{+\infty}T(t,\cdot)(\mathrm{Id}-P(\cdot))B(\cdot)U(\cdot,s),~~~t\geq s.
\label{UTP}
\end{split}
\end{eqnarray}
\label{zmlm1}
\end{lm}

{\bf Proof.} Clearly, if the function $U(t,s)_{t\geq s}$ satisfies \eqref{UTP},
then it is a solution of \eqref{TxABx}. We must demonstrate that the operator $L$ defined by
\begin{eqnarray*}
(LU)(t,s)=\!\!\!\!\!\!\!\!&&T(t,s)P(s)+{\cal I}^\alpha_{s}T(t,\cdot)P(\cdot)B(\cdot)U(\cdot,s)\\
&&+{\cal I}^\alpha_{+\infty}T(t,\cdot)(\mathrm{Id}-P(\cdot))B(\cdot)U(\cdot,s),~~~t\geq s,
\end{eqnarray*}
has a unique fixed point in $W$. It follows from \eqref{TPK} that
\begin{eqnarray*}
\|(LU)(t,s)\|\leq\!\!\!\!\!\!\!\!&& \|T(t,s)P(s)\|+{\cal I}^\alpha_{s}\|T(t,\cdot)P(\cdot)\|\|B(\cdot)\|\|U(\cdot,s)\|
\\
&&-{\cal I}^\alpha_{+\infty}\|T(t,\cdot)(\mathrm{Id}-P(\cdot))\|\|B(\cdot)\|\|U(\cdot,s)\|
\\
\leq\!\!\!\!\!\!\!\!&& \hat{N}_1\frac{E_\alpha(\hat{\beta}_1,s)}{E_\alpha(\hat{\beta}_1,t)}E_\alpha(\epsilon,s)
+\delta\left(\frac{\hat{N}_1}{\hat{\beta}_1}+\frac{\hat{N}_2}{\hat{\beta}_2}\right) \|U\|_\alpha E_\alpha(\epsilon,s).
\end{eqnarray*}
Combining \eqref{uslts} with \eqref{the}, we obtain
\begin{eqnarray*}
\|LU\|_\alpha\leq \hat{N}_1 +\theta\|U\|_\alpha<\infty,
\end{eqnarray*}
this implies that the operator $L: W\to W$ is well defined.
Analogously to the computation above, we have
\begin{eqnarray*}
\|LU_1-LU_2\|_\alpha\leq \theta\|U_1-U_2\|_\alpha,~~~U_1,U_2\in W,
\end{eqnarray*}
which shows that $L$ is a contraction since $\theta<1$.
Then there exists a unique $U\in W$ such that $LU=U$.
Therefore, Lemma \ref{zmlm1} is proved.
$\qquad\Box$

Now we explain that the bounded solutions exhibit the following property.

\begin{lm}
For each $t\geq \tau \geq s$ in $\mathbb{R}_+$,
\begin{eqnarray*}
U(t,\tau)U(\tau,s)=U(t,s).
\end{eqnarray*}
\label{zmlm2}
\end{lm}

{\bf Proof.} From \eqref{UTP} and \eqref{TPPT}, for some $\tau\in \mathbb{R}_+$ we can calculate that
\begin{eqnarray*}
U(t,\tau)U(\tau,s)=\!\!\!\!\!\!\!\!&&T(t,s)P(s)+{\cal I}^\alpha_{s}T(t,\tau)P(\tau)B(\tau)U(\tau,s)
\\
&&+{\cal I}^\alpha_{\tau}T(t,\cdot)P(\cdot)B(\cdot)U(\cdot,\tau)U(\tau,s)
\\
&&+{\cal I}^\alpha_{+\infty}T(t,\cdot)(\mathrm{Id}-P(\cdot))B(\cdot)U(\cdot,\tau)U(\tau,s),~~~t\geq \tau \geq s.
\end{eqnarray*}
Let $H(t,\tau):=U(t,\tau)U(\tau,s)-U(t,s)$ for $t\geq \tau \geq s$, this yields
\begin{eqnarray}
H(t,\tau)={\cal I}^\alpha_{\tau}T(t,\cdot)P(\cdot)B(\cdot)H(\cdot,s)
+{\cal I}^\alpha_{+\infty}T(t,\cdot)(\mathrm{Id}-P(\cdot))B(\cdot)H(\cdot,s).
\label{HtI}
\end{eqnarray}
Define operator $\mathcal{K}$ as
\begin{eqnarray*}
(\mathcal{K}\hat{H})(t,\tau):={\cal I}^\alpha_{\tau}T(t,\cdot)P(\cdot)B(\cdot)\hat{H}(\cdot,s)
+{\cal I}^\alpha_{+\infty}T(t,\cdot)(\mathrm{Id}-P(\cdot))B(\cdot)\hat{H}(\cdot,s),
\end{eqnarray*}
for any $\hat{H}\in W$ and $t\geq \tau$.
It follows from the identity above and \eqref{TPK} that
\begin{eqnarray*}
\|(\mathcal{K}\hat{H})(t,\tau))\|\leq \!\!\!\!\!\!\!\!&&{\cal I}^\alpha_{\tau}\|T(t,\cdot)P(\cdot)\|\|B(\cdot)\|\|\hat{H}(\cdot,s)\|
\\
&&-{\cal I}^\alpha_{+\infty}\|T(t,\cdot)(\mathrm{Id}-P(\cdot))\|\|B(\cdot)\|\|\hat{H}(\cdot,s)\|\\
\leq \!\!\!\!\!\!\!\!&& \delta\left(\frac{\hat{N}_1}{\hat{\beta}_1}+\frac{\hat{N}_2}{\hat{\beta}_2}\right)\|\hat{H}\|_\alpha E_\alpha(\epsilon,s).
\end{eqnarray*}
By \eqref{uslts}, we have
\begin{eqnarray*}
\|\mathcal{K}\hat{H}\|_\alpha\leq \theta\|\hat{H}\|_\alpha<\infty,
\end{eqnarray*}
then $\mathcal{K}: W\to W$ is well defined for $t\geq \tau$.
Similarly to the calculation above, we attain
\begin{eqnarray*}
\|\mathcal{K}\hat{H}_1-\mathcal{K}\hat{H}_2\|_\alpha\leq \theta\|\hat{H}_1-\hat{H}_2\|_\alpha,~~~\hat{H}_1,\hat{H}_2\in W.
\end{eqnarray*}
Because of hypothesis \eqref{the}, $\mathcal{K}$ is a contraction.
Thus, there is a unique $\hat{H}\in W$ such that $\mathcal{K}\hat{H}=\hat{H}$.
On the other hand, we know that $0\in W$ satisfies \eqref{HtI} and $\mathcal{K}0=0$.
By Lemma \ref{zmlm1}, we assert $H=\hat{H}=0$ for $t\geq \tau \geq s$ in $\mathbb{R}_+$.
Therefore, Lemma \ref{zmlm2} is proved.
$\qquad\Box$

Step 2: Establishment of projections $\hat{P}(t)$ in \eqref{Tts}.
Given constant $\iota\in \mathbb{R}_+$, for any $t\geq \iota$ in $\mathbb{R}_+$, we consider the following linear operator
\begin{eqnarray}
\hat{P}(t):=\hat{T}(t,\iota)U(\iota,\iota)\hat{T}(\iota,t),
\label{PTU}
\end{eqnarray}
where $\hat{T}(t,s)$ is the evolution operator associated to \eqref{TxABx}.
Clearly, the operator $\hat{P}(t)$ may depend on $\iota$,
and $U(\iota,\iota)U(\iota,\iota)=U(\iota,\iota)$ by Lemma \ref{zmlm2}.
The following lemma illustrates the commutativity of projections $\hat{P}(t)$ as formula \eqref{Tts}.
\begin{lm}
For any $t\in \mathbb{R}_+$, the operator $\hat{P}(t)$ is a projection satisfying \eqref{Tts}.
\label{zmlm3}
\end{lm}

{\bf Proof.}
By the details above and {\bf (F1)-(F2)}, we derive
\begin{eqnarray*}
\hat{P}(t)\hat{P}(t)=\!\!\!\!\!\!\!\!&&\hat{T}(t,\iota)U(\iota,\iota)\hat{T}(\iota,t)\hat{T}(t,\iota)U(\iota,\iota)\hat{T}(\iota,t)
\\
=\!\!\!\!\!\!\!\!&&\hat{T}(t,\iota)U(\iota,\iota)U(\iota,\iota)\hat{T}(\iota,t)=\hat{P}(t),
\end{eqnarray*}
then $\hat{P}(t)$ is a projection. Furthermore, for $t\geq s$ we can calculate that
\begin{eqnarray*}
\hat{T}(t,s)\hat{P}(s)=\!\!\!\!\!\!\!\!&& \hat{T}(t,s)\hat{T}(s,\iota)U(\iota,\iota)\hat{T}(\iota,t)\hat{T}(t,s)=\hat{P}(t)\hat{T}(t,s).
\end{eqnarray*}
This completes the proof of Lemma \ref{zmlm3}.
$\qquad\Box$

Step 3: Characterization of bounded solutions.
The following two lemmas propose the nonuniform projection integral equation
and its property respectively.

\begin{lm}
For some $s\in \mathbb{R}_+$, if $z\in C_b([s,+\infty), Z)$ is a solution of \eqref{TxABx} with $z(s)=z_s$, then
\begin{eqnarray*}
z(t)=\!\!\!\!\!\!\!\!&&T(t,s)P(s)z_s+{\cal I}^\alpha_{s}T(t,\cdot)P(\cdot)B(\cdot)z(\cdot)
+{\cal I}^\alpha_{+\infty}T(t,\cdot)(\mathrm{Id}-P(\cdot))B(\cdot)z(\cdot).
\end{eqnarray*}
\label{zmlm4}
\end{lm}

The proof of this lemma is similar to the method of Lemma \ref{efds} when $\epsilon< \min\{\hat{\beta}_1,\hat{\beta}_2\}$ holds.

\begin{lm}
For some $s\in \mathbb{R}_+$, if the function $\hat{P}(\cdot)\hat{T}(\cdot,s)\in C_b([s,+\infty) , \mathcal{B}(Z))$, then
\begin{eqnarray}
\begin{split}
\hat{P}(t)\hat{T}(t,s)=&T(t,s)P(s)\hat{P}(s)+{\cal I}^\alpha_{s}T(t,\cdot)P(\cdot)B(\cdot)\hat{P}(\cdot)\hat{T}(\cdot,s)
\\
&+{\cal I}^\alpha_{+\infty}T(t,\cdot)(\mathrm{Id}-P(\cdot))B(\cdot)\hat{P}(\cdot)\hat{T}(\cdot,s).
\label{PtTts}
\end{split}
\end{eqnarray}
\label{zmlm5}
\end{lm}

{\bf Proof.} For a given $\iota\in \mathbb{R}_+$, it follows from Lemma \ref{zmlm1} that
the function $U(t,\iota)\xi$ is a solution of \eqref{TxABx}
with initial value $U(\iota,\iota)\xi$ for any $\xi\in Z$ and $t\geq \iota$.
By \eqref{PTU} and \eqref{Tts}, we gain $U(t,\iota)=\hat{T}(t,\iota)U(\iota,\iota)$, and
\begin{eqnarray*}
\hat{P}(t)\hat{T}(t,s)=\!\!\!\!\!\!\!\!&&\hat{T}(t,s)\hat{P}(s)=\hat{T}(t,s)\hat{T}(s,\iota)U(\iota,\iota)\hat{T}(\iota,s)
\\
=\!\!\!\!\!\!\!\!&&\hat{T}(t,\iota)U(\iota,\iota)\hat{T}(\iota,s)=U(t,\iota)\hat{T}(\iota,s).
\label{PTT}
\end{eqnarray*}
Thus, the equation \eqref{TxABx} has solution in the form of $U(t,\iota)\xi$ as follows
\begin{eqnarray*}
z(t)=\hat{P}(t)\hat{T}(t,s)\xi=U(t,\iota)\hat{T}(\iota,s)\xi,~~~\xi\in Z.
\end{eqnarray*}
Observing that the above solution is bounded for $t\geq s$, and
\begin{eqnarray*}
z(s)=U(s,\iota)\hat{T}(\iota,s)\xi=\hat{P}(s)\hat{T}(s,s)\xi=\hat{P}(s)\xi,
\end{eqnarray*}
we employ Lemma \ref{zmlm4} to complete the proof of Lemma \ref{zmlm5}.
$\qquad\Box$

The following Lemma is the projected integral inequality in the case of nonuniform Mittag-Leffler dichotomy,
and the method of its proof can be referred to the Lemma \ref{efbds1} and Corollary \ref{efbds2}.

\begin{lm}
Given $s\in \mathbb{R}_+$. Assume that the functions $u\in C_b([s,+\infty),\mathbb{R}_+)$ and $v\in C_b([0,s],\mathbb{R}_+)$
respectively satisfy the following inequalities
\begin{eqnarray}
\begin{split}
u(t)\leq& \hat{N}_1\frac{E_\alpha(\hat{\beta}_1,s)}{E_\alpha(\hat{\beta}_1,t)}E_\alpha(\epsilon,s)u_s
+\frac{\delta \hat{N}_1}{E_\alpha(\hat{\beta}_1,t)}{\cal I}^\alpha_{s}E_\alpha(\hat{\beta}_1,t)u(t)
\\
&-\delta \hat{N}_2E_\alpha(\hat{\beta}_2,t){\cal I}^\alpha_{+\infty}\frac{u(t)}{E_\alpha(\hat{\beta}_2,t)},
~~~t\geq s\geq0,
\label{un}
\end{split}
\\
\begin{split}
v(t)\leq& \hat{N}_2\frac{E_\alpha(\hat{\beta}_2,t)}{E_\alpha(\hat{\beta}_2,s)}E_\alpha(\epsilon,s)v_s
+\frac{\delta \hat{N}_1}{E_\alpha(\hat{\beta}_1,t)}{\cal I}^\alpha_{0}E_\alpha(\hat{\beta}_1,t)v(t)
\\
&-\delta \hat{N}_2 E_\alpha(\hat{\beta}_2,t){\cal I}^\alpha_{s}\frac{v(t)}{E_\alpha(\hat{\beta}_2,t)},
~~~s\geq t\geq0,
\label{vn}
\end{split}
\end{eqnarray}
where $u_s:=u(s)$ and $v_s:=v(s)$. If
\begin{eqnarray*}
\theta:=\delta\Big{(}\frac{\hat{N}_1}{\hat{\beta}_1}+\frac{\hat{N}_2}{\hat{\beta}_2}\Big{)}<1,
\end{eqnarray*}
then there exist positive constants $K_i$ and $\lambda_i(i=1,2)$ such that
\begin{eqnarray*}
u(t)\leq\!\!\!\!\!\!\!\!\!&& K_1\frac{E_\alpha(\lambda_1,s)}{E_\alpha(\lambda_1,t)}E_\alpha(\epsilon,s)u_s,~~~t\geq s\geq0,
\\
v(t)\leq\!\!\!\!\!\!\!\!\!&& K_2\frac{E_\alpha(\lambda_2,t)}{E_\alpha(\lambda_2,s)}E_\alpha(\epsilon,s)v_s,~~~s\geq t\geq0,
\end{eqnarray*}
where $K_i=\cfrac{\hat{N}_i}{1-\theta}$, $\lambda_i=\hat{\beta}_i-\cfrac{\delta \hat{N}_i}{1-\theta}$.
\label{zmlm6}
\end{lm}

Step 4: Norm bounds of evolution operator.
We verify that the norms of the operators $\hat{T}(t,s)|\mathrm{Im}\hat{P}(s)$
and $\hat{T}(t,s)|\mathrm{Im}(\mathrm{Id}-\hat{P}(s))$ are bounded.
\begin{lm}
For any $t\geq s$ in $\mathbb{R}_+$, the first inequality in \eqref{TPKE} holds.
\label{zmlm7}
\end{lm}

{\bf Proof.} Given $\xi\in Z$, and for $t\geq s\geq 0$, assume that
\begin{eqnarray*}
u(t):=\|\hat{P}(t)\hat{T}(t,s)\xi\|,
\end{eqnarray*}
then $u_s=\|\hat{P}(s)\xi\|$.
By Lemma \ref{zmlm5}, we know that $u(t)$ is bounded and satisfies \eqref{un}.
It follows from Lemma \ref{zmlm6} that
\begin{eqnarray*}
\|\hat{P}(t)\hat{T}(t,s)\xi\|\leq K_1\frac{E_\alpha(\lambda_1,s)}{E_\alpha(\lambda_1,t)}E_\alpha(\epsilon,s)\|\hat{P}(s)\xi\|,
~~~t\geq s\geq 0,
\end{eqnarray*}
where $K_1$ and $\lambda_1$ are given in Lemma \ref{zmlm6}.
Again by Lemma \ref{zmlm3}, we gain
\begin{eqnarray*}
\hat{P}(t)\hat{T}(t,s)=\hat{T}(t,s)\hat{P}(s)=\hat{T}(t,s)\hat{P}(s)\hat{P}(s).
\end{eqnarray*}
Taking $\mu:=\hat{P}(s)\xi$, it yields that
\begin{eqnarray*}
\|\hat{T}(t,s)\hat{P}(s)\mu\|\leq K_1\frac{E_\alpha(\lambda_1,s)}{E_\alpha(\lambda_1,t)}E_\alpha(\epsilon,s)\|\mu\|,
~~~t\geq s\geq 0.
\end{eqnarray*}
Therefore, we can obtain the desired inequality.
$\qquad\Box$

\begin{lm}
For any $s\geq t$ in $\mathbb{R}_+$, the second inequality in \eqref{TPKE} holds.
\label{zmlm8}
\end{lm}

{\bf Proof.} By analogy with Lemma \ref{zmlm5},
we need to attain an equation for $(\mathrm{Id}-\hat{P}(t))\hat{T}(t,s)$ via Lemma \ref{zmlm3}.
Actually, from the variation of constants formula \eqref{csbygs}, we have
\begin{eqnarray*}
\hat{T}(t,s)=T(t,s)+{\cal I}^\alpha_{s}T(t,\cdot)B(\cdot)\hat{T}(\cdot,s).
\end{eqnarray*}
Let function $w(t):=\hat{T}(t,\iota)(\mathrm{Id}-\hat{P}(\iota))$ for some $\iota\in \mathbb{R}_+$, then
\begin{eqnarray}
w(t)=T(t,\iota)(\mathrm{Id}-\hat{P}(\iota))+{\cal I}^\alpha_{\iota}T(t,\cdot)B(\cdot)w(\cdot).
\label{legs1}
\end{eqnarray}
From \eqref{UTP} and \eqref{PTU} with $t=s=\iota$, we calculate that
\begin{eqnarray*}
\hat{P}(\iota)=U(\iota,\iota)=P(\iota)+{\cal I}^\alpha_{+\infty}T(\iota,\cdot)(\mathrm{Id}-P(\cdot))B(\cdot)U(\cdot,\iota).
\end{eqnarray*}
Pre-projecting $P(\iota)$ on both hands sides of the above identity,
we acquire $P(\iota)\hat{P}(\iota)=P(\iota)$, and
\begin{eqnarray}
(\mathrm{Id}-P(\iota))(\mathrm{Id}-\hat{P}(\iota))=\mathrm{Id}-\hat{P}(\iota).
\label{PPP}
\end{eqnarray}
Combining \eqref{legs1} with \eqref{PPP}, and replacing $t$ with $s$, we derive
\begin{eqnarray*}
T(t,s)(\mathrm{Id}-P(s))w(s)=\!\!\!\!\!\!\!\!&&T(t,\iota)(\mathrm{Id}-P(\iota))(\mathrm{Id}-\hat{P}(\iota))
\nonumber\\
&&+{\cal I}^\alpha_{\iota}T(t,s)(\mathrm{Id}-P(s))B(s)w(s)
\nonumber\\
=\!\!\!\!\!\!\!\!&&T(t,\iota)(\mathrm{Id}-\hat{P}(\iota))+{\cal I}^\alpha_{\iota}T(t,s)(\mathrm{Id}-P(s))B(s)w(s).
\label{legs2}
\end{eqnarray*}
It follows from \eqref{legs1} and the identity above that
\begin{eqnarray}
w(t)=\!\!\!\!\!\!\!\!&&T(t,s)(\mathrm{Id}-P(s))w(s)+{\cal I}^\alpha_{\iota}T(t,\cdot)B(\cdot)w(\cdot)
\nonumber\\
&&-{\cal I}^\alpha_{\iota}T(t,s)(\mathrm{Id}-P(s))B(s)w(s)
\nonumber\\
=\!\!\!\!\!\!\!\!&&T(t,s)(\mathrm{Id}-P(s))w(s)+{\cal I}^\alpha_{\iota}T(t,\cdot)P(\cdot)B(\cdot)w(\cdot)
\nonumber\\
&&+{\cal I}^\alpha_{s}T(t,\cdot)(\mathrm{Id}-P(\cdot))B(\cdot)w(\cdot).
\label{legs3}
\end{eqnarray}
On the other hand, by Lemma \ref{zmlm3}, we attain
\begin{eqnarray}
(\mathrm{Id}-\hat{P}(t))\hat{T}(t,s)=\hat{T}(t,s)(\mathrm{Id}-\hat{P}(s)).
\label{PTts}
\end{eqnarray}
Recalling the function $w(\tau)$,
we get $w(\tau)\hat{T}(\iota,s)=(\mathrm{Id}-\hat{P}(\tau))\hat{T}(\tau,s)$.
Post-multiplying $\hat{T}(\iota,s)$ on both hands sides of \eqref{legs3},
this implies
\begin{eqnarray}
\begin{split}
(\mathrm{Id}-\hat{P}(t))\hat{T}(t,s)=&T(t,s)(\mathrm{Id}-P(s))(\mathrm{Id}-\hat{P}(s))
\\
&+{\cal I}^\alpha_{\iota}T(t,\cdot)P(\cdot)B(\cdot)(\mathrm{Id}-\hat{P}(\cdot))\hat{T}(\cdot,s)
\\
&+{\cal I}^\alpha_{s}T(t,\cdot)(\mathrm{Id}-P(\cdot))B(\cdot)(\mathrm{Id}-\hat{P}(\cdot))\hat{T}(\cdot,s).
\label{IDPt}
\end{split}
\end{eqnarray}
Fixed $\xi\in Z$, we consider $v(t):=\|\hat{T}(t,s)(\mathrm{Id}-\hat{P}(s))\xi\|$ for $s\geq t\geq 0$
and $v_s=\|(\mathrm{Id}-\hat{P}(s))\xi\|$.
According to \eqref{legs3} and \eqref{PTts},
it is well known that the function $v(t)$ satisfies the inequality \eqref{vn}.
Employing Lemma \ref{zmlm6} and the similar proof to Lemma \ref{zmlm7},
we easily acquire desired inequality and complete the proof.
$\qquad\Box$

In conclusion, Lemmas \ref{zmlm3}, \ref{zmlm7} and \ref{zmlm8} all derive Theorem \ref{thmfyz} together.
$\qquad\Box$

The following Lemma will help to prove Theorem \ref{thmfyz1}.

\begin{lm}
For any $t\in \mathbb{R}_+$, if constant $\delta$ described as in Theorem \ref{thmfyz1} is small enough, then
\begin{eqnarray}
\|\hat{P}(t)\|\leq 4\hat{N}E_\alpha(\epsilon, t),~~~\|\mathrm{Id}-\hat{P}(t)\|\leq 4\hat{N}E_\alpha(\epsilon, t).
\label{PtNE}
\end{eqnarray}
\label{zmlm9}
\end{lm}

{\bf Proof.} Replacing $s$ by $t$ and pre-multiplying $(\mathrm{Id}-P(t))$ on both hands sides of \eqref{PtTts}, we have
\begin{eqnarray}
(\mathrm{Id}-P(t))\hat{P}(t)={\cal I}^\alpha_{+\infty}T(t,\cdot)(\mathrm{Id}-P(\cdot))B(\cdot)\hat{P}(\cdot)\hat{T}(\cdot,t).
\label{prgs1}
\end{eqnarray}
It follows from Lemmas \ref{zmlm7} and \ref{zmlm3} that for $\tau\geq t\geq 0$,
\begin{eqnarray}
\|\hat{P}(\tau)\hat{T}(\tau,t)\|=\|\hat{T}(\tau,t)\hat{P}(t)\hat{P}(t)\|
\leq K_1\frac{E_\alpha(\lambda_1,t)}{E_\alpha(\lambda_1,\tau)}E_\alpha(\epsilon,t)\|\hat{P}(t)\|.
\label{prgs2}
\end{eqnarray}
By \eqref{prgs1} and \eqref{TPK} we calculate that
\begin{eqnarray}
\|(\mathrm{Id}-P(t))\hat{P}(t)\|\leq\!\!\!\!\!\!\!\!&&-{\cal I}^\alpha_{+\infty}\|T(t,\cdot)(\mathrm{Id}-P(\cdot))\|\|B(\cdot)\|\|\hat{P}(\cdot)\hat{T}(\cdot,t)\|
\nonumber\\
\leq\!\!\!\!\!\!\!\!&&  -E_\alpha(\hat{\beta}_2+\lambda_1+\epsilon,t)\|\hat{P}(t)\|
{\cal I}^\alpha_{+\infty}\frac{\delta K_1 \hat{N}_2}{E_\alpha(\hat{\beta}_2+\lambda_1+\epsilon,t)}
\nonumber\\
\leq\!\!\!\!\!\!\!\!&& \frac{\delta K_1 \hat{N}_2}{\hat{\beta}_2+\lambda_1-\epsilon}\|\hat{P}(t)\|,
\label{prgs3}
\end{eqnarray}
where constant $\epsilon$ was chosen as satisfying $\epsilon<\min\{\hat{\beta}_1,\hat{\beta}_2\}$
in order to guarantee the above denominator $\hat{\beta}_2+\lambda_1-\epsilon>0$.
Analogously to \eqref{prgs1}, replacing $t$ with $s$ and pre-multiplying $P(t)$ on both hands sides of \eqref{IDPt},
we attain
\begin{eqnarray}
P(t)(\mathrm{Id}-\hat{P}(t))={\cal I}^\alpha_{\iota}T(t,\cdot)P(\cdot)B(\cdot)(\mathrm{Id}-\hat{P}(\cdot))\hat{T}(\cdot,t).
\label{prgs4}
\end{eqnarray}
Using Lemma \ref{zmlm8}, for $t\geq \tau \geq 0$ this implies
\begin{eqnarray}
\|(\mathrm{Id}-\hat{P}(\tau))\hat{T}(\tau,t)\|\leq
K_2\frac{E_\alpha(\lambda_2,\tau)}{E_\alpha(\lambda_2,t)}E_\alpha(\epsilon,t)\|\mathrm{Id}-\hat{P}(t)\|.
\label{prgs5}
\end{eqnarray}
From \eqref{prgs4} and \eqref{TPK} one can compute that
\begin{eqnarray}
\|P(t)(\mathrm{Id}-\hat{P}(t))\|\leq\!\!\!\!\!\!\!\!&&{\cal I}^\alpha_{\iota}\|T(t,\cdot)P(\cdot)\|\|B(\cdot)\|\|(\mathrm{Id}-\hat{P}(\cdot))\hat{T}(\cdot,t)\|
\nonumber\\
\leq\!\!\!\!\!\!\!\!&& \frac{\delta K_2 \hat{N}_1}{E_\alpha(\hat{\beta}_1+\lambda_2-\epsilon,t)}\|\mathrm{Id}-\hat{P}(t)\|
{\cal I}^\alpha_{\iota}E_\alpha(\hat{\beta}_1+\lambda_2-\epsilon,t)
\nonumber\\
\leq\!\!\!\!\!\!\!\!&& \frac{\delta K_2 \hat{N}_1}{\hat{\beta}_1+\lambda_2-\epsilon}\|\mathrm{Id}-\hat{P}(t)\|,
\label{prgs6}
\end{eqnarray}
where the chosen constant $\epsilon<\min\{\hat{\beta}_1,\hat{\beta}_2\}$ similarly.
Obviously,
\begin{eqnarray*}
\hat{P}(t)-P(t)=\!\!\!\!\!\!\!\!&&(\mathrm{Id}-P(t))\hat{P}(t)-P(t)(\mathrm{Id}-\hat{P}(t)).
\end{eqnarray*}
Taking $\hat{N}:=\max\{\hat{N}_1,\hat{N}_2\}$ and $\hat{\beta}:=\min\{\hat{\beta}_1,\hat{\beta}_2\}$
and combining \eqref{prgs3} with \eqref{prgs6}, we gain
\begin{eqnarray}
\|\hat{P}(t)-P(t)\|\leq\!\!\!\!\!\!\!\!&& \frac{\delta K_1 \hat{N}_2}{\hat{\beta}_2+\lambda_1-\epsilon}\|\hat{P}(t)\|
+\frac{\delta K_2 \hat{N}_1}{\hat{\beta}_1+\lambda_2-\epsilon}\|\mathrm{Id}-\hat{P}(t)\|
\nonumber\\
\leq\!\!\!\!\!\!\!\!&& \hat{\eta}(\|\hat{P}(t)\|+\|\mathrm{Id}-\hat{P}(t)\|),
\label{prgs7}
\end{eqnarray}
where
$$\hat{\eta}=\cfrac{\delta\hat{N}^2\hat{\beta}}{2\hat{\beta}^2-5\delta\hat{N}\hat{\beta}-\epsilon(\hat{\beta}-2\delta\hat{N})}.$$
Moreover, by \eqref{TPK} with $t=s$, it is easy to obtain that
\begin{eqnarray*}
\|P(t)\|\leq \hat{N}E_\alpha(\epsilon,t),~~~\|Q(t)\|\leq \hat{N}E_\alpha(\epsilon, t).
\end{eqnarray*}
From \eqref{prgs7}, this yields
\begin{eqnarray*}
\|\hat{P}(t)\|\leq\!\!\!\!\!\!\!\!&& \|\hat{P}(t)-P(t)\|+\|P(t)\|
\\
\leq\!\!\!\!\!\!\!\!&& \hat{\eta}(\|\hat{P}(t)\|+\|\mathrm{Id}-\hat{P}(t)\|)
+\hat{N}E_\alpha(\epsilon,t).
\end{eqnarray*}
Since $\|(\mathrm{Id}-\hat{P}(t))-(\mathrm{Id}-P(t))\|=\|\hat{P}(t)-P(t)\|$, we also derive
\begin{eqnarray*}
\|(\mathrm{Id}-\hat{P}(t))\|\leq\!\!\!\!\!\!\!\!&& \|\hat{P}(t)-P(t)\|+\|\mathrm{Id}-P(t)\|
\\
\leq\!\!\!\!\!\!\!\!&& \hat{\eta}(\|\hat{P}(t)\|+\|\mathrm{Id}-\hat{P}(t)\|)
+\hat{N}E_\alpha(\epsilon,t).
\end{eqnarray*}
They together imply that
\begin{eqnarray*}
\|\hat{P}(t)\|+\|\mathrm{Id}-\hat{P}(t)\|\leq 2\hat{\eta}(\|\hat{P}(t)\|+\|\mathrm{Id}-\hat{P}(t)\|)
+2\hat{N}E_\alpha(\epsilon,t),
\end{eqnarray*}
and
\begin{eqnarray*}
\|\hat{P}(t)\|+\|\mathrm{Id}-\hat{P}(t)\|\leq \frac{2\hat{N}E_\alpha(\epsilon,t)}{1-2\hat{\eta}}.
\end{eqnarray*}
Choose $\hat{\eta}<1/4$, then
\begin{eqnarray*}
\|\hat{P}(t)\|+\|\mathrm{Id}-\hat{P}(t)\|\leq 4\hat{N}E_\alpha(\epsilon,t),
\end{eqnarray*}
yielding Lemma \ref{zmlm9}.
$\qquad\Box$

Finally, we end this paper with the proof of roughness for nonuniform Mittag-Leffler dichotomy.

{\bf Proof of Theorem \ref{thmfyz1}.}
From \eqref{prgs2} and \eqref{PtNE}, we show that
\begin{eqnarray*}
\|\hat{P}(\tau)\hat{T}(\tau,t)\|\leq\!\!\!\!\!\!\!\!&&
\frac{\hat{N}\hat{\beta}}{\hat{\beta}-2\delta\hat{N}}\frac{E_\alpha(\hat{\lambda},t)}{E_\alpha(\hat{\lambda},\tau)}E_\alpha(\epsilon,t)\|\hat{P}(t)\|
\\
\leq\!\!\!\!\!\!\!\!&& \frac{4\hat{N}^2\hat{\beta}}{\hat{\beta}-2\delta\hat{N}}
\frac{E_\alpha(\hat{\lambda},t)}{E_\alpha(\hat{\lambda},\tau)}E_\alpha(2\epsilon,t),~~~\tau\geq t\geq 0,
\end{eqnarray*}
where $\hat{\lambda}=\hat{\beta}-\cfrac{\delta\hat{N}\hat{\beta}}{\hat{\beta}-2\delta\hat{N}}$.
Analogously, it follows from \eqref{prgs5} and \eqref{PtNE} that
\begin{eqnarray*}
\|(\mathrm{Id}-\hat{P}(\tau))\hat{T}(\tau,t)\|\leq\!\!\!\!\!\!\!\!&&
\frac{4\hat{N}^2\hat{\beta}}{\hat{\beta}-2\delta\hat{N}}
\frac{E_\alpha(\hat{\lambda},\tau)}{E_\alpha(\hat{\lambda},t)}E_\alpha(2\epsilon,t),~~~t\geq \tau\geq 0.
\end{eqnarray*}
Therefore, we can acquire the desired inequalities like \eqref{TPK},
and the proof is completed.
$\qquad\Box$


\end{document}